\documentclass[aos,nameyear,dvips]{arximspdf}

\doi{10.1214/10-AOS807}
\volume{38}
\issue{5}
\pubyear{2010}
\firstpage{2823}
\lastpage{2856}

\makeatletter

\newtheorem{theorem}{Theorem}[section]
\newtheorem{lemma}{Lemma}[section]
\newproclaim{df}{Definition}[section]

\makeatother

\begin{document}
\begin{frontmatter}

\title{Trajectory averaging for stochastic approximation MCMC algorithms\thanksref{T1}}
\runtitle{Trajectory averaging for SAMCMC}

\thankstext{T1}{Supported in part by NSF Grants DMS-06-07755, CMMI-0926803 and the Award KUS-C1-016-04
made by King Abdullah University of
Science and Technology (KAUST).}

\begin{aug}
\author[A]{\fnms{Faming} \snm{Liang}\corref{}\ead[label=e1]{fliang@stat.tamu.edu}}
\runauthor{F. Liang}
\affiliation{Texas A\&M University}
\address[A]{Department of Statistics\\
Texas A\&M University\\
College Station, Texas 77843-3143\\
USA\\
\printead{e1}} 
\end{aug}

\received{\smonth{5} \syear{2009}}
\revised{\smonth{11} \syear{2009}}

%
\begin{abstract}
The subject of stochastic approximation was founded by Robbins and Monro [\textit{Ann. Math. Statist.} \textbf{22} (\citeyear{RM51}) 400--407].
After five decades of continual development, it has
developed into an important area in systems control and optimization,
and it has also served as a prototype for the development of adaptive
algorithms for on-line estimation and control of stochastic systems.
Recently, it has been used in statistics with Markov chain Monte Carlo
for solving maximum likelihood estimation problems and for general
simulation and optimizations. In this paper, we first show that the
trajectory averaging estimator is asymptotically efficient for the
stochastic approximation MCMC (SAMCMC) algorithm under mild conditions,
and then apply this result to the stochastic approximation Monte Carlo
algorithm [Liang, Liu and Carroll \textit{J. Amer. Statist. Assoc.} \textbf{102} (\citeyear{LLC07}) 305--320].
The application of the trajectory averaging estimator to other
stochastic approximation MCMC algorithms, for example, a stochastic
approximation MLE algorithm for missing data problems, is also
considered in the paper.
\end{abstract}

%

\end{frontmatter}

\section{Introduction}

Robbins and Monro (\citeyear{RM51}) introduced the stochastic approximation
algorithm to solve the
integration equation
%
%
\begin{equation} \label{problemeq1}
h(\theta)=\int_{{\mathcal X}} H(\theta,x) f_{\theta}(x) \,dx =0,
\end{equation}
where $\theta\in\Theta\subset\mathbb{R}^{d_{\theta}}$ is a
parameter vector
and $f_{\theta}(x)$, $x \in{\mathcal X}\subset\mathbb{R}^{d_x}$, is
a density
function depending on $\theta$.
The $d_{\theta}$ and $d_x$ denote the dimensions of $\theta$ and $x$,
respectively.
The stochastic approximation algorithm is an iterative recursive
algorithm, whose each iteration
consists of two steps:

\subsection*{Stochastic approximation algorithm}
\begin{itemize}
\item Generate $X_{k+1} \sim f_{\theta_k}(x)$, where $k$ indexes the iteration.
\item Set $\theta_{k+1}=\theta_k+ a_k H(\theta_k, X_{k+1})$, where
$a_k>0$ is called the gain factor.
\end{itemize}

The stochastic approximation algorithm is often studied by rewriting
it as follows:
%
%
\begin{equation} \label{moteq1}
\theta_{k+1}=\theta_k+a_k[ h(\theta_k)+\varepsilon_{k+1}],
\end{equation}
where $h(\theta_k)=\int_{{\mathcal X}} H(\theta_k,x) f_{\theta_k}(x)\,dx$
corresponds to the mean effect
of $H(\theta_k,\break X_{k+1})$, and $\varepsilon_{k+1}=H(\theta_k$,
$X_{k+1})-h(\theta_k)$ is called the
observation noise. In the literature of stochastic approximation,
$h(\theta)$ is also called the
mean field function.
It is well known that the optimal convergence rate of (\ref{moteq1})
can be
achieved with $a_k=-F^{-1}/k$, where $F=\partial h(\theta^*)/\partial
\theta$, and
$\theta^*$ denotes the zero point of $h(\theta)$.
In this case, (\ref{moteq1}) is reduced to Newton's algorithm.
Unfortunately, it is often
impossible to use this algorithm, as the matrix $F$ is generally unknown.

Although an optimal convergence rate of $\theta_k$ cannot be obtained
in general,
in a sequence of fundamental papers \citet{R88}, \citet{P90}
and \citet{PJ92} showed that the trajectory averaging
estimator is
asymptotically efficient; that is,
$\bar{\theta}_n=\sum_{k=1}^n \theta_k/n$
can converge in distribution to a normal random variable with
mean $\theta^*$ and covariance matrix $\Sigma$,
where $\Sigma$ is the smallest possible covariance matrix in an appropriate
sense. The trajectory averaging estimator requires $\{a_k\}$ to be
relatively large,
decreasing slower than $O(1/k)$.
As discussed by \citet{PJ92}, trajectory averaging
is based on a paradoxical principle:
a slow algorithm having less than optimal convergence rate must be averaged.

Recently, the trajectory averaging technique has been further explored
in a variety of papers
[see, e.g., \citet{C93}, Kushner and Yang (\citeyear{KY93}, \citeyear{KY95}), \citet{DR97},
Wang, Chong and Kulkarni (\citeyear{WCK97}), \citet{TLEC99},
\citet{P00} and \citet{KY03}] with different assumptions
for the observation noise.
However, up to our knowledge, it has not yet
been explored for stochastic approximation MCMC (SAMCMC) algorithms
[\citet{BMP90}, \citet{C02}, \citet{KY03}, \citet{AMP05},
\citet{AM06}]. The stochastic approximation MCMC
algorithms refer to a class of
stochastic approximation algorithms for which the sample is generated
at each iteration
via a Markov transition kernel; that is, $\{x_{k+1}\}$ is generated
via a family of
Markov transition kernel $\{ P_{\theta_k}(x_k,\cdot) \}$ controlled by
$\{\theta_k \}$.
Recently, the stochastic approximation MCMC algorithms have been used
in statistics
for solving maximum likelihood estimation problems [Younes (\citeyear{Y89}, \citeyear{Y99}),
\citet{MB91}, \citet{GK98}, \citet{GZ01}],
and for general simulation and optimizations [\citet{LLC07}, \citet{AL10}].
It is worth to point out that in comparison with conventional MCMC
algorithms, for example,
the Metropolis--Hastings algorithm [\citet{Metal53}, \citet{H70}],
parallel tempering [\citet{G91}], and simulated tempering [\citet{MP92},
\citet{GT95}],
the stochastic approximation Monte Carlo (SAMC) algorithm [\citet{LLC07}]
has significant advantages in simulations of complex systems for which
the energy landscape
is rugged. As explained later (in Section \ref{sec3}), SAMC is essentially
immune to the local trap problem due to
its self-adaptive nature inherited from the stochastic approximation algorithm.
SAMC has been successfully applied to many statistical problems, such as
$p$-value evaluation for resampling-based tests [\citet{YL09}],
Bayesian model selection [\citet{L09}, \citet{AL10}]
and spatial model estimation [Liang (\citeyear{L07a})], among others.

In this paper, we explore the theory of trajectory averaging for stochastic
approximation MCMC algorithms, motivated by their wide applications.
Although Chen (\citeyear{C93}, \citeyear{C02}) considered the case where the observation
noise can
be state dependent, that is, the observation noise $\varepsilon_{k+1}$
depends on
$\theta_0, \ldots, \theta_k$, their results are not directly
applicable to the
stochastic approximation MCMC algorithms due to some reasons as
explained in Section \ref{sec5}.
The theory established by \citet{KY03}
can potentially be extended to the stochastic approximation MCMC
algorithm, but,
as mentioned in Kushner and Yin [(\citeyear{KY03}), page 375] the extension is not
straightforward
and more work needs to be done to deal with the complicated structure
of the Markov
transition kernel.
In this paper, we propose a novel decomposition of the observation
noise for the
stochastic approximation MCMC algorithms. Based on the proposed
decomposition, we show
the trajectory averaging estimator is asymptotically
efficient for the stochastic approximation MCMC algorithms,
and then apply this result
to the SAMC algorithm.
These results are presented in Lemma \ref{lem3}, Theorems
\ref{efftheorem} and
\ref{samcavetheorem}, respectively.
The application of the trajectory averaging technique to other stochastic
approximation MCMC algorithms, for example, a stochastic approximation
MLE algorithm for
missing data problems, is also considered in the paper.

The remainder of this paper is organized as follows. In Section \ref{sec2}, we
present our main theoretical result that the trajectory averaging
estimator is
asymptotically efficient for the stochastic approximation MCMC algorithms.
In Section \ref{sec3}, we apply the trajectory averaging technique to the SAMC algorithm.
In Section \ref{sec4}, we apply the trajectory averaging technique to a
stochastic approximation MLE
algorithm for missing data problems.
In Section \ref{sec5}, we conclude the paper with a brief discussion.

\eject
\section{Trajectory averaging for a general stochastic
approximation MCMC algorithm}\label{sec2}

\subsection{A varying truncation stochastic approximation MCMC algorithm}\label{sec21}

To show the convergence of the stochastic approximation algorithm,
restrictive conditions on the observation noise and mean field function
are required. For example, one often assumes the noise to be mutually
independent or to be a martingale difference sequence, and imposes
a sever restriction on the growth rate of the mean field function.
These conditions are usually not satisfied in practice.
See Chen [(\citeyear{C02}), Chapter 1] for more discussions on this issue.
To remove the growth rate restriction on the mean field function and to
weaken the conditions imposed on noise, \citet{CZ86} proposed a
varying truncation version for the stochastic approximation algorithm.
The convergence of the modified algorithm can be shown for a wide
class of the mean filed
function under a truly weak condition on noise; see, for example, \citet{CGG88}
and \citet{AMP05}. The latter gives a proof for the convergence
of the modified algorithm with Markov state-dependent noise
under some conditions that are easy to verify.

Following \citet{AMP05}, we consider the
following varying truncation stochastic approximation MCMC algorithm.
Let $\{\mathcal{K}_s, s \geq0\}$ be a sequence of compact subsets of
$\Theta$
such that
%
%
\begin{equation} \label{truncationseteq}
\bigcup_{s \geq0} \mathcal{K}_s=\Theta\quad \mbox{and}\quad \mathcal{K}_s
\subset\mbox{int}(\mathcal{K}_{s+1}),\qquad s \geq0,
\end{equation}
where int($A$) denotes the interior of set $A$. Let $\{a_k\}$
and $\{b_k\}$ be two monotone, nonincreasing, positive sequences.
Let ${\mathcal X}_0$ be a subset of ${\mathcal X}$, and let $\mathcal
{T}\dvtx {\mathcal X}\times
\Theta\rightarrow{\mathcal X}_0 \times\mathcal{K}_0$ be a
measurable function which
maps a point $(x,\theta)$ in ${\mathcal X}\times\Theta$ to a random
point in
${\mathcal X}_0 \times\mathcal{K}_0$; that is,
both $x$ and $\theta$ will be reinitialized in ${\mathcal X}_0 \times
\mathcal{K}_0$.
As shown in Lemma \ref{lem3}, for the stochastic approximation MCMC
algorithm, when the number of iterations becomes large, the
observation noise $ \varepsilon_k $ can be
decomposed as
%
%
\begin{equation} \label{rev.eq1}
\varepsilon_{k}=e_k+\nu_k+\varsigma_k,
\end{equation}
where $\{e_k\}$ forms a martingale difference sequence, and the
expectation of the other
two terms will go to zero in certain forms.
In Theorems \ref{contheorem} and \ref{efftheorem}, we show that $\{
e_k\}$
leads to the asymptotic normality of the trajectory averaging estimator
$\bar{\theta}_k$, and $\{\nu_k\}$ and $\{ \varsigma_k \}$
can vanish or be ignored when the asymptotic distribution of $\bar
{\theta}_k$ is considered.

Let $\sigma_k$ denote the number of truncations performed until
iteration $k$ and $\sigma_0=0$.
The varying truncation stochastic approximation MCMC algorithm starts
with a random choice of
$(\theta_0,x_0)$ in the space $\mathcal{K}_0 \times{\mathcal X}_0$,
and then iterates
between the following steps:

\subsubsection*{Varying truncation stochastic approximation MCMC algorithm}
\begin{itemize}
\item Draw sample $x_{k+1}$ with a Markov transition kernel, $P_{\theta
_k}$, which admits $f_{\theta_k}(x)$
as the invariant distribution.

\item Set $\theta_{k+{1/2}}=\theta_k+a_k H(\theta_k,x_{k+1})$.

\item If $\|\theta_{k+{1/2}} -\theta_k \| \leq b_k$ and
$\theta_{k+{1/2}} \in\mathcal{K}_{\sigma_k}$, where $\|z\|$
denote the
Euclidean norm of the vector $z$,
then set $(\theta_{k+1},x_{k+1})=(\theta_{k+{1/2}},x_{k+1})$ and
$\sigma_{k+1}=\sigma_k$; otherwise, set $(\theta
_{k+1},x_{k+1})=\mathcal{T}
(\theta_{k},x_k)$ and
$\sigma_{k+1}=\sigma_k+1$.
\end{itemize}

As depicted by the algorithm, the varying truncation mechanism works
in an adaptive manner
as follows: when the current estimate of the parameter wanders outside
the active
truncation set or when the difference between two successive estimates
is greater
than a time-dependent threshold, then the algorithm is reinitialized
with a smaller
initial value of the gain factor and a larger truncation set. This mechanism
enables the algorithm to select an appropriate gain factor sequence
and an appropriate
starting point, and thus to
confine the recursion to a compact set; that is, the number of reinitializations
is almost surely finite for every $(\theta_0,x_0) \in\mathcal{K}_0
\times
{\mathcal X}_0$.
This result is formally stated in Theorem \ref{lem50}, which
plays a crucial role for establishing asymptotic efficiency of
the trajectory averaging estimator.

Regarding the varying truncation scheme, one
can naturally propose many variations. For example,
one may not change the truncation set when only the
condition $\|\theta_{k+{1/2}} -\theta_k \| \leq b_k$ is violated,
and, instead of jumping forward in a unique gain factor sequence,
one may start with a different gain factor sequence (smaller than the
previous one)
when the reinitialization occurs.
In either case, the proof for the theorems presented in Section
\ref{SAMCMCtheory}
follows similarly.

\subsection{Theoretical results on the trajectory averaging estimator}
\label{SAMCMCtheory}

The asymptotic efficiency of $\bar{\theta}_k$ can be analyzed under
the following conditions.

\subsubsection*{Lyapunov condition on $h(\theta)$}

Let $\langle x, y \rangle$ denote the Euclidean inner product.

\begin{enumerate}[(A$_1$)]
\item[(A$_1$)] $\Theta$ is an open set, the function $h\dvtx \Theta
\rightarrow\mathbb{R}^d$ is
continuous, and there exists a continuously differentiable function
$v\dvtx \Theta\rightarrow
[0,\infty)$ such that:
\begin{longlist}
\item There exists $M_0>0$ such that
%
%
\begin{equation} \label{solutionseteq}
\mathcal{L}=\{ \theta\in\Theta, \langle\nabla v(\theta),h(\theta
)\rangle=0
\} \subset
\{ \theta\in\Theta, v(\theta) <M_0 \}.
\end{equation}

\item There exists $M_1 \in(M_0, \infty)$ such that $\mathcal
{V}_{M_1}$
is a compact set,
where $\mathcal{V}_{M}=\{ \theta\in\Theta, v(\theta) \leq M \}$.

\item For any $\theta\in\Theta\setminus\mathcal{L}$,
$\langle
\nabla
v(\theta),h(\theta)\rangle< 0$.

\item The closure of $v(\mathcal{L})$ has an empty interior.
\end{longlist}
\end{enumerate}

This condition assumes the existence of a global Lyapunov function $v$
for the mean field $h$.
If $h$ is a gradient field, that is, $h=-\nabla J$ for some lower
bounded real-valued and
differentiable function $J(\theta)$, then $v$ can be set to $J$,
provided that $J$ is
continuously differentiable. This is typical for stochastic
optimization problems, for example, machine
learning [\citet{T97}], where a continuously
differentiable objective function
$J(\theta)$ is minimized.

\subsubsection*{Stability condition on $h(\theta)$}

\begin{enumerate}[(A$_2$)]
\item[(A$_2$)] The mean field function $h(\theta)$ is measurable and
locally bounded.
There exist a stable matrix $F$ (i.e., all eigenvalues of $F$ are with
negative real parts),
$\gamma>0$, $\rho\in(0,1]$, and a constant $c$ such that, for any
$\theta^* \in\mathcal{L}$,
\[
\| h(\theta)-F(\theta-\theta^*)\| \leq c
\|\theta-\theta^*\|^{1+\rho}\qquad \forall\theta\in\{ \theta\dvtx
\|\theta-\theta^*\| \leq\gamma\},
\]
where $\mathcal{L}$ is defined in (\ref{solutionseteq}).
\end{enumerate}

This condition constrains the behavior of the mean field function around
the solution points. It makes the trajectory averaging estimator sensible
both theoretically and practically.
If $h(\theta)$ is differentiable, the matrix $F$ can be chosen to be
the partial derivative of $h(\theta)$, that is, $\partial h(\theta
)/\partial\theta$.
Otherwise, certain approximation may be needed.

\subsubsection*{Drift condition on the transition kernel $P_{\theta}$}
Before giving details of this condition, we first define some terms
and notation.
Assume that a transition kernel $P_{\theta}$ is irreducible,
aperiodic, and has a stationary distribution on
a sample space denoted by ${\mathcal X}$.
A set $\mathbf{C}\subset{\mathcal X}$ is said to be small if there
exist a
probability measure $\nu$ on ${\mathcal X}$,
a positive integer $l$ and $\delta>0$ such that
\[
P_{\theta}^l(x,A) \geq\delta\nu(A)\qquad \forall x \in\mathbf{C},
\forall A
\in\mathcal{B}_{{\mathcal X}},
\]
where $\mathcal{B}_{{\mathcal X}}$ is the Borel set defined on
${\mathcal X}$.
A function $V\dvtx {\mathcal X}\rightarrow[1,\infty)$ is said to be a drift
function outside $\mathbf{C}$ if there exist
positive constants $\lambda<1$ and $b$ such that
\[
P_{\theta} V(x) \leq\lambda V(x)+b I(x \in\mathbf{C})\qquad \forall x
\in{\mathcal X},
\]
where $P_{\theta} V(x)=\int_{{\mathcal X}} P_{\theta}(x,y) V(y)\,d y$.
For a function $g\dvtx {\mathcal X}\rightarrow\mathbb{R}^d$, define the norm
\[
\| g\|_V=\sup_{x \in{\mathcal X}} \frac{\|g(x)\|}{V(x)}
\]
and define the set $\mathcal{L}_V=\{g\dvtx {\mathcal X}\rightarrow\mathbb
{R}^d, {\sup_{x \in
{\mathcal X}}}
\|g\|_V < \infty\}$.
Given the terms and notation introduced above, the drift condition can
be specified as follows.

\begin{enumerate}[(A$_3$)]
\item[(A$_3$)] For any given $\theta\in\Theta$,
the transition kernel $P_{\theta}$ is irreducible and aperiodic.
In addition, there exists a
function $V\dvtx {\mathcal X}\rightarrow[1,\infty)$ and a constant
$\alpha\geq
2$ such that
for any compact subset $\mathcal{K}\subset\Theta$:
%
\begin{longlist}
\item There exist a set $\mathbf{C}\subset{\mathcal X}$,
an integer $l$, constants $0<\lambda<1$, $b$, $\varsigma$, $\delta>0$
and a probability
measure $\nu$ such that
%
%
\begin{eqnarray}
\label{Drieq11}
 \sup_{\theta\in\mathcal{K}}P_{\theta}^l V^{\alpha}(x)
&\leq&
\lambda V^{\alpha}(x)+b I(x \in\mathbf{C})\qquad
\forall x \in{\mathcal X}, \\
\label{Drieq12}
 \sup_{\theta\in\mathcal{K}} P_{\theta} V^{\alpha
}(x)&\leq&
\varsigma V^{\alpha}(x)\qquad
\forall x \in{\mathcal X},\\
\label{Drieq13}
 \inf_{\theta\in\mathcal{K}} P_{\theta}^l(x, A) &\geq&
\delta\nu(A)\qquad
\forall x \in\mathbf{C}, \forall A \in\mathcal{B}_{{\mathcal X}}.
\end{eqnarray}

\item There exists a constant $c>0$ such that, for all $x \in
{\mathcal X}$,
%
%
\begin{eqnarray}
\label{Drieq21}
 {\sup_{\theta\in\mathcal{K}}} \|H(\theta,x)\|_V &\leq& c,
\\
\label{Direq22}
 {\sup_{(\theta,\theta')\in\mathcal{K}}} \|H(\theta
,x)-H(\theta
',x)\|_V
&\leq& c \|\theta-\theta'\|.
\end{eqnarray}
\item There exists a constant $c>0$ such that, for all $(\theta
,\theta')\in\mathcal{K}\times\mathcal{K}$,
%
%
\begin{eqnarray} \label{Drieq31}
   \| P_{\theta}g-P_{\theta'} g\|_V &\leq& c \|g \|_V
\|\theta-\theta' \|\qquad \forall g \in\mathcal{L}_V, \\
\label{Drieq32}
  \| P_{\theta}g-P_{\theta'} g\|_{V^{\alpha}} &\leq& c \|
g\|
_{V^{\alpha}}
\|\theta-\theta' \|\qquad \forall g \in\mathcal{L}_{V^{\alpha}}.
\end{eqnarray}

\end{longlist}
\end{enumerate}

Assumption (A$_3$)(i) is classical in the literature of Markov
chain. It implies the existence
of a stationary distribution $f_{\theta}(x)$ for all $\theta\in
\Theta
$ and $V^{\alpha}$-uniform
ergodicity [\citet{AMP05}].
Assumption (A$_3$)(ii) gives conditions on the bound of $H(\theta,x)$.
This is a critical condition for the observation noise.
As seen later in Lemmas \ref{lem1} and \ref{lem3}, it directly leads
to the boundedness
of some terms decomposed from the observation noise.
For some algorithms, for example, SAMC, for which $H(\theta,x)$ is a
bounded function, the drift function
can be simply set as $V(x)=1$.

\subsubsection*{Conditions on the step-sizes}

\begin{enumerate}[(A$_4$)]
\item[(A$_4$)] The sequences $\{a_k\}$ and $\{b_k\}$ are
nonincreasing, positive and satisfy the conditions:
%
%
\begin{eqnarray} \label{coneq1}
\sum_{k=1}^{\infty} a_k&=&\infty,\qquad
\lim_{k \rightarrow\infty} (k a_k) =\infty,\nonumber\\[-8pt]\\[-8pt]
\frac{a_{k+1}-a_k}{a_k}&=&o(a_{k+1}),\qquad
b_k=O\bigl(a_k^{({1+\tau})/{2}}\bigr),\nonumber
\end{eqnarray}
for some $\tau\in(0,1]$,
%
%
\begin{equation} \label{coneq2}
\sum_{k=1}^{\infty} \frac{a_k^{(1+\tau)/2}}{\sqrt{k}}< \infty,
\end{equation}
and for some constants $\alpha\geq2$
as defined in condition (A$_3$),
%
%
\begin{equation} \label{coneq0003}
\sum_{i=1}^{\infty} \{ a_i b_i +(b_i^{-1} a_i)^{\alpha}
\} < \infty.
\end{equation}
\end{enumerate}

It follows from (\ref{coneq2}) that
\[
\sum_{i=[{k/2}]}^k \frac{a_i^{(1+\tau)/2}}{\sqrt{i}} =o(1),
\]
where $[z]$ denotes the integer part of $z$. Since $a_k$ is
nonincreasing, we have
\[
a_k^{(1+\tau)/2} \sum_{i=[{k}/{2}]}^k \frac{1}{\sqrt{i}} =o(1),
\]
and thus $a_k^{(1+\tau)/2} \sqrt{k}=o(1)$, or $a_k=O(k^{-\eta})$ for
$\eta\in
(\frac{1}{2},1)$.
For instance, $a_k=C_1/k^{\eta}$ for some constants $C_1>0$ and $\eta
\in
(\frac{1}{2},1)$, then we can set $b_k=C_2/k^{\xi}$ for some constants
$C_2>0$ and
$\xi\in(\frac{1}{2}, \eta-\frac{1}{\alpha})$, which satisfies
(\ref
{coneq1}) and (\ref{coneq0003}).
Under this setting, the existence
of $\tau$ is obvious.

Theorem \ref{lem50} concerns the convergence of the general
stochastic approximation MCMC algorithm. The proof follows directly from
Theorems 5.4, 5.5 and Proposition 6.1 of \citet{AMP05}.
\begin{theorem} \label{lem50}
Assume  conditions \textup{(A$_1$), (A$_3$)} and \textup{(A$_4$)} hold.
Let $k_{\sigma}$ denote the iteration number at which
the $\sigma$th truncation occurs in the stochastic approximation MCMC
simulation.
Let ${\mathcal X}_0 \subset{\mathcal X}$ be such that $\sup_{x\in
{\mathcal X}_0} V(x)<\infty$
and that
$\mathcal{K}_0 \subset\mathcal{V}_{M_0}$, where $\mathcal{V}_{M_0}$
is defined in \textup{(A$_1$)}.
Then there exists almost surely a number, denoted by $\sigma_s$, such that
$k_{\sigma_s}<\infty$ and $k_{\sigma_s+1}=\infty$; that is, $\{
\theta_k\}$
has no truncation for $k \geq k_{\sigma_s}$, or mathematically,
\[
\theta_{k+1}=\theta_k+a_k H(\theta_k,x_{k+1})\qquad \forall k \geq
k_{\sigma_s}.
\]
In addition, we have
\[
\theta_k \to\theta^*\qquad\mbox{a.s.}
\]
for some point $\theta^* \in\mathcal{L}$.
\end{theorem}

Theorem \ref{contheorem} concerns the asymptotic normality
of $\bar{\theta}_k$.
\begin{theorem} \label{contheorem}
Assume  conditions \textup{(A$_1$), (A$_2$), (A$_3$)} and \textup{(A$_4$)} hold.
Let ${\mathcal X}_0 \subset{\mathcal X}$ be such that $\sup_{x\in
{\mathcal X}_0} V(x)<\infty$
and that
$\mathcal{K}_0 \subset\mathcal{V}_{M_0}$, where $\mathcal{V}_{M_0}$
is defined in \textup{(A$_1$)}.
Then
\[
\sqrt{k} (\bar{\theta}_k-\theta^*) \longrightarrow N(\mathbf{0},\Gamma)
\]
for some point $\theta^* \in\Theta$,
where $\Gamma=F^{-1}Q (F^{-1})^T$, $F=\partial h(\theta^*)/\partial
\theta$ is
negative definite,
$Q=\lim_{k \rightarrow\infty} E (e_k e_k^T)$, and $e_k$ is as
defined in
(\ref{rev.eq1}).
\end{theorem}

Below we consider the asymptotic efficiency of $\bar{\theta}_k$.
As already mentioned, the asymptotic efficiency of the
trajectory averaging estimator has been studied by quite a few authors.
\citet{TLEC99} gives the following definition for the
asymptotic efficient
estimator that can be resulted from a stochastic approximation algorithm.
\begin{df} \label{effdef}
Consider the stochastic approximation algorithm (\ref{moteq1}).
Let $\{Z_n\}_{n \geq0}$, given as a function of $\{\theta_n\}_{n \geq0}$,
be a sequence of estimators of $\theta^*$.
The algorithm $\{Z_n\}_{n \geq0}$ is said to be asymptotically
efficient if
%
%
\begin{equation} \label{dfeq2}
\sqrt{n}(Z_n-\theta^*) \longrightarrow N ( \mathbf{0}, F^{-1}\tilde{Q}
(F^{-1})^T ),
\end{equation}
where $F=\partial h(y^*)/\partial y$, and $\tilde{Q}$ is the
asymptotic covariance matrix
of $(1/\sqrt{n})\times\break \sum_{k=1}^n \varepsilon_k$.
\end{df}

As mentioned in \citet{TLEC99},
$\tilde{Q}$ is the smallest possible limit covariance
matrix that an estimator based on the stochastic approximation
algorithm (\ref{moteq1})
can achieve. If $\theta_k \rightarrow\theta^*$ and $\{\varepsilon_k\}$
forms or asymptotically
forms a martingale difference sequence, then we have $\tilde
{Q}=\lim_{k \rightarrow\infty}
E (\varepsilon_k \varepsilon_k^T)$. In the next theorem, we show that
the asymptotic covariance matrix $Q$ established in Theorem~\ref{contheorem} is the same as
$\tilde{Q}$, and thus the trajectory averaging estimator $\bar
{\theta}_k$ is asymptotically efficient.
%
\begin{theorem} \label{efftheorem}
Assume  conditions \textup{(A$_1$), (A$_2$), (A$_3$)} and \textup{(A$_4$)} hold.
Let ${\mathcal X}_0 \subset{\mathcal X}$ be such that $\sup_{x\in
{\mathcal X}_0} V(x)<\infty$
and that
$\mathcal{K}_0 \subset\mathcal{V}_{M_0}$, where $\mathcal{V}_{M_0}$
is defined in \textup{(A$_1$)}. Then
$\bar{\theta}_k$ is asymptotically efficient.
\end{theorem}

As implied by Theorem \ref{efftheorem}, the convergence rate of $\bar
{\theta}_k$,
which is measured by the asymptotic covariance matrix $\Gamma$,
is independent of the choice of the gain factor sequence as long as
the condition
(A$_4$) is satisfied.
The asymptotic efficiency of $\bar{\theta}_k$ can also be interpreted
in terms of
Fisher information theory.
Refer to Pelletier [(\citeyear{P00}), Section 3] and the references therein for more
discussions on this issue.

Trajectory averaging enables smoothing of the behavior of the
algorithm but at the
same time, it slows down the numerical convergence because it takes
longer for
the algorithm to forget the first iterates. An alternative idea would
be to consider
moving window averaging algorithms,
see, for example, \citet{KY93} and Kushner and Yin
(\citeyear{KY03}), Chapter 11. Extension of their results to the general
stochastic approximation MCMC algorithm will be of great interest.

\eject
\section{Trajectory averaging for the stochastic
approximation Monte Carlo algorithm}\label{sec3}

\subsection{The SAMC algorithm}\label{sec31}

Suppose that we are interested in sampling from the following distribution
%
%
\begin{equation} \label{tareq1}
f(x)=c \psi(x),\qquad x \in{\mathcal X},
\end{equation}
where $c$ is an unknown constant, ${\mathcal X}\subset\mathbb
{R}^{d_x}$ is the
sample space.
The basic idea of SAMC stems from the Wang--Landau algorithm [\citet{WL01},
\citet{L05}]
and can be briefly explained as follows.
Let $E_1, \ldots, E_m$ denote a partition of ${\mathcal X}$, and let
$\omega
_i=\int_{E_i} \psi(x)\,dx$
for $i=1, \ldots, m$.
SAMC seeks to draw sample from the trial distribution
%
%
\begin{equation} \label{eq11}
f_{\omega}(x) \propto\sum_{i=1}^m \frac{\pi_i \psi(x)}{\omega_i}
I_{\{
x \in E_i\}},
\end{equation}
where $\pi_i$'s are prespecified constants such that $\pi_i>0$ for
all $i$
and \mbox{$\sum_{i=1}^m \pi_i=1$}, and
$I_{\{x \in E_i\}}=1$ if $x \in E_i$ and 0 otherwise.
For example,
if the sample space is partitioned according to the energy function
into the following subregions:
$E_1=\{x\dvtx -\log(\psi(x))<u_1\}$, $E_2=\{x\dvtx u_1 \leq-\log(\psi(x))
<u_2\}, \ldots,
E_m=\break\{x\dvtx -\log(\psi(x))> u_{m-1}\}$, where $-\infty< u_1 < \cdots<
u_{m-1}< \infty$ are the user-specified numbers,
then sampling from $f_{\omega}(x)$ would result in a random walk (by
viewing each subregion as a ``point'')
in the space of energy with
each subregion being sampled with probability $\pi_i$. Here, without
loss of generality,
we assume that each subregion is
unempty; that is, assuming $\int_{E_i} \psi(x) \,dx >0$ for all
$i=1,\ldots, m$.
Therefore, sampling from (\ref{eq11}) essentially avoids the
local-trap problem suffered by
the conventional MCMC algorithms. This is attractive, but
$\omega_i$'s are unknown.
SAMC provides a dynamic way to estimate $\omega_i$'s under the framework
of the stochastic approximation MCMC algorithm.

In what follows we describe how $\omega$ can be estimated by SAMC.
Since $f_{\omega}(x)$ is invariant with respect to a scale change
of $\omega$, it suffices to
estimate $\omega_1,\ldots,\omega_{m-1}$ by fixing $\omega_m$ to a
known constant provided $\omega_m>0$.
Let $\theta_{k}^{(i)}$ denote the working estimate of $\log(\omega
_i/\pi_i)$ obtained at
iteration $k$, and let $\theta_k=(\theta_{k}^{(1)}, \ldots,\theta
_{k}^{(m-1)})$. Why this
reparameterization is used will be explained at the end of this subsection.
Let $\{a_k\}$ denote the gain factor sequence, and let
$\{\mathcal{K}_s, s \geq0\}$ denote a sequence of compact
subsets of $\Theta$ as defined in (\ref{truncationseteq}).
For this algorithm, $\{ \mathcal{K}_s, s \geq0\}$ can be chosen as
follows. Define
%
%
\begin{equation} \label{vfunctioneq}
v(\theta)=-\log\Biggl( 1- \frac{1}{2} \sum_{j=1}^{m-1} \biggl(\frac
{S_j}{S}-\pi
_j\biggr)^2 \Biggr),
\end{equation}
where $S_i=\int_{E_i} \psi(x) \,dx/\exp(\theta^{(i)})$ for
$i=1,\ldots,
m-1$, and
$S=\sum_{i=1}^{m-1} S_i+\int_{E_i} \psi(x) \,dx$. Clearly, $v(\theta)$
is continuous in $\theta$,
and $\mathcal{V}_M=\{\theta\dvtx v(\theta) \leq M\}$ for any $M \in(0,
\infty)$
forms a compact subset of $\Theta$.
Therefore, $\{ \mathcal{V}_{M_s}, s \geq0\}$, $0<M_0 < M_1 < \cdots
$, is an
appropriate choice
of $\{\mathcal{K}_s, s \geq0\}$.
For the SAMC algorithm, as seen below, $\| H(\theta_k,X_{k+1})\|= \|
(I_{\{x_{k+1} \in E_1\}} -\pi_1,
\ldots, I_{\{x_{k+1} \in E_{m-1}\}} -\pi_{m-1})^T \|$ is bounded by
the constant $\sqrt{2}$, so
we can set the drift function $V(x)=1$. Hence, the initial sample
$x_0$ can be drawn arbitrarily from ${\mathcal X}_0={\mathcal X}$,
while leaving the condition $\sup_{x\in{\mathcal X}_0} V(x) <\infty$ holds.
In summary, SAMC starts with an initial estimate of $\theta_0 \in
\mathcal{K}_0$,
and a random sample drawn arbitrarily from the
space ${\mathcal X}$, and then iterates between the following steps.

\textit{SAMC algorithm}.
\begin{enumerate}[(a)]
\item[(a)](Sampling.) Simulate a sample $x_{k+1}$ by a single MH update
with the target distribution
%
%
\begin{equation} \label{tardiseq}
f_{\theta_k}(x) \propto\sum_{i=1}^{m-1} \frac{\psi(x)}{e^{\theta
_{k}^{(i)}}} I_{\{x \in E_i\}}
+ \psi(x)I_{\{x \in E_m\}},
\end{equation}
provided that $E_m$ is nonempty. In practice, $E_m$ can be replaced by
any other unempty subregion.
\begin{enumerate}[(a.2)]
\item[(a.1)] Generate $y$ according to a proposal distribution $q(x_k,y)$.
\item[(a.2)] Calculate the ratio
\[
r=e^{\theta_k^{(J(x_k))} -\theta_k^{(J(y))}} \frac{\psi(y) q(y,
x_k)}{\psi(x_k) q(x_k,y)},
\]
where $J(z)$ denotes the index of the subregion that the sample $z$
belongs to.
\item[(a.3)] Accept the proposal with probability $\min(1,r)$. If it is
accepted, set $x_{k+1}=y$;
otherwise, set $x_{k+1}=x_k$.
\end{enumerate}

\item[(b)] (Weight updating.) Set
%
%
\begin{equation} \label{weiupdate}
\theta_{k+{1}/{2}}^{(i)}=
\theta_k^{(i)}+a_{k+1} \bigl( I_{\{x_{k+1} \in E_i\}}-\pi_i \bigr),\qquad i=1,\ldots,m-1.
\end{equation}

\item[(c)] (Varying truncation.) If $\theta_{k+{1}/{2}}
\in\mathcal{K}_{\sigma_k}$, then set
$(\theta_{k+1},x_{k+1})=(\theta_{k+{1}/{2}},x_{k+1})$ and
$\sigma_{k+1}=\sigma_k$; otherwise, set
$(\theta_{k+1},x_{k+1})=\mathcal{T}(\theta_{k},x_k)$ and
$\sigma_{k+1}=\sigma_k+1$, where $\sigma_k$ and $\mathcal{T}(\cdot
,\cdot )$ are as defined in Section \ref{sec2}.
\end{enumerate}

SAMC sampling is driven by its self-adjusting mechanism, which, consequently,
implies the superiority of SAMC in sample space exploration.
The self-adjusting mechanism can be explained as follows:
if a subregion is visited at \mbox{iteration $k$}, $\theta_k$ will be updated
accordingly
such that the probability that this subregion (other subregions) will be
revisited at the next iterations will decrease (increase).
Mathematically,\vspace*{-2pt} if $x_{k+1} \in E_i$,
then $\theta_{k+{1}/{2}}^{(i)} \leftarrow\theta_{k}^{(i)}+
a_{k+1} (1-\pi_i)$
and $\theta_{k+{1}/{2}}^{(j)} \leftarrow\theta_{k}^{(j)}-a_{k+1}
\pi_j$ for $j \ne i$.
Note that the linear adjustment on $\theta$ transforms to a
multiplying adjustment
on $\omega$. This also explains why SAMC works on the logarithm of
$\omega$.
Working on the logarithm
enables $\omega$ to be adjusted quickly according to the distribution
of the samples.
Otherwise, learning of $\omega$ would be very slow due to the linear
nature of
stochastic approximation. Including $\pi_i$ in the transformation
$\log(\omega_i/\pi_i)$ facilitates our computation, for example, the
ratio $r$
in step (a.2).

The self-adjusting mechanism has led to successful applications of
SAMC for many hard computational problems, including phylogenetic tree
reconstruction [Cheon and Liang (\citeyear{CL07}, \citeyear{CL09})], neural network training
[Liang (\citeyear{L07b})], Bayesian network learning [\citet{LZ09}],
among others.

\subsection{Trajectory averaging for SAMC}\label{sec32}

To show that the trajectory averaging estimator is asymptotically efficient
for SAMC, we assume the following conditions.

\begin{enumerate}[(C$_1$)]
\item[(C$_1$)] The MH transition kernel used in the sampling step
satisfies the drift condition (A$_3$).
\end{enumerate}

To ensure the drift condition to be satisfied, \citet{LLC07}
restrict
the sample space ${\mathcal X}$ to be a compact set, assume $f(x)$ to be
bounded away from
0 and $\infty$, and choose the proposal distribution
$q(x,y)$ to satisfy the local positive condition:
for every $x \in{\mathcal X}$, there exist positive
$\varepsilon_1$ and $\varepsilon_2$ such that
%
%
\begin{equation} \label{proposalcons}
\|x-y\| \leq\varepsilon_1 \quad\Longrightarrow \quad q(x,y) \geq\varepsilon_2.
\end{equation}
If the compactness condition on ${\mathcal X}$ is removed,
we may need to impose some constraints on the tails of the target distribution
$f(x)$ and the proposal distribution $q(x,y)$ as done by \citet{AMP05}.

\begin{enumerate}[(C$_2$)]
\item[(C$_2$)] The sequence $\{a_k\}$ satisfies the following conditions:
\begin{eqnarray*}
\sum_{k=1}^{\infty} a_k&=&\infty,\qquad \lim_{k \to\infty} (k
a_k)=\infty,\\
\frac{a_{k+1}-a_k}{a_k}&=&o(a_{k+1}),\qquad
\sum_{k=1}^{\infty} \frac
{a_k^{(1+\tau)/2}}{\sqrt{k}} < \infty
\end{eqnarray*}
for some $\tau\in(0,1]$.
\end{enumerate}
For the SAMC algorithm, as previously discussed, $\| H(\theta
_k,X_{k+1})\|$
is bounded by the constant $\sqrt{2}$, so
we can set $V(x)=1$ and set $\alpha$ to any a large number in
condition (A$_3$).
Furthermore, given a choice of $a_k=O(k^{-\eta})$ for some $\eta\in(1/2,1)$,
there always exists a sequence $\{b_k\}$, for example,
$b_k=2 a_k^{(1+\tau)/2}$ for some $\tau\in(0,1]$,
such that the inequality $\| \theta_{k+{1}/{2}} -\theta_k \|=\|
a_k H(\theta_k,X_{k+1})\| \leq b_k$ holds for all iterations.
Hence, a specification of the sequence $\{b_k\}$ can be omitted for
the SAMC algorithm.

Theorem \ref{samcconvergence} concerns the convergence of SAMC.
In the first part, it states that $k_{\sigma_s}$ is almost surely
finite; that is, $\{\theta_k\}$ can be included in a compact set
almost surely. In the second part, it states the
convergence of $\theta_k$ to a solution $\theta^*$. We note that for SAMC,
the same convergence result has been established by \citet{LLC07}
under (C$_1$) and a relaxed condition of (C$_2$),
where $\{a_k\}$ is allowed to decrease at a rate of $O(1/k)$.
Since the focus of this paper is on the asymptotic efficiency of $\bar
{\theta}_k$,
the convergence of $\{\theta_k\}$ is only stated under a slower
decreasing rate of $\{a_k\}$.
We also note that for SAMC, we have assumed, without loss of
generality, that all
subregions are unempty.
For the empty subregions,
no adaptation of $\{\theta_k\}$ occurs for the corresponding
components in the run.
Therefore, the convergence of $\{\theta_k\}$ should only be measured
for the components corresponding to the nonempty subregions.
\begin{theorem} \label{samcconvergence}
Assume  conditions \textup{(C$_1$)} and \textup{(C$_2$)} hold.
Then there exists (a.s.) a number, denoted by $\sigma_s$, such that
$k_{\sigma_s}<\infty$, $k_{\sigma_s+1}=\infty$, and $\{\theta_k\}$ given
by the SAMC algorithm has no truncation for $k \geq k_{\sigma_s}$,
that is,
%
%
\begin{equation} \label{convergeq0}
\theta_{k+1}=\theta_k+a_k H(\theta_k,x_{k+1})\qquad \forall k \geq
k_{\sigma_s}
\end{equation}
and
%
%
\begin{equation} \label{convergeq}
\theta_k \rightarrow\theta^*\qquad\mbox{a.s.},
\end{equation}
where $H(\theta_k,x_{k+1})= (I_{\{x_{k+1}\in E_1\}}-\pi_1, \ldots,
I_{\{x_{k+1}\in E_{m-1}\}}-\pi_{m-1} )^T$, and
$\theta^*= ( \log(\omega_1/\pi_1)-\log(\omega_m/\pi_m), \ldots,
\log(\omega_{m-1}/\pi_{m-1})-\log(\omega_m/\pi_m) )^T$.
\end{theorem}

Theorem \ref{samcavetheorem} concerns the asymptotic normality and
efficiency of
$\bar{\theta}_k$.
\begin{theorem} \label{samcavetheorem}
Assume  conditions \textup{(C$_1$)} and \textup{(C$_2$)}. Then $\bar{\theta}_k$
is asymptotically efficient; that is,
\[
\sqrt{k} (\bar{\theta}_k-\theta^*) \longrightarrow N(\mathbf{0},\Gamma)
\qquad\mbox{as $k \rightarrow\infty$,}
\]
where $\Gamma=F^{-1}Q (F^{-1})^T$, $F=\partial h(\theta^*)/\partial
\theta$ is negative definite and $Q=\break\lim_{k \rightarrow\infty} E(e_k e_k^T)$.
\end{theorem}

The above theorems address some theoretical issues of SAMC. For
practical issues,
please refer to \citet{LLC07}, where issues, such as
how to partition the sample space, how to choose the desired sampling
distribution, and how to diagnose the convergence, have
been discussed at length. An issue particularly related to the trajectory
averaging estimator is the length of the burn-in period. To remove
the effect of the early iterates, the following estimator:
\[
\bar{\theta}_k^{(b)}=\frac{1}{k-k_0} \sum_{i=k_0+1}^k \theta_i,
\]
instead of $\bar{\theta}_k$, is often used in practice, where $k_0$ is
the so-called
length of the burn-in period. It is obvious that the choice of $k_0$
should be based on the diagnosis for the convergence of the simulation.
Just like monitoring
convergence of MCMC simulations, monitoring convergence of SAMC simulations
should be based on multiple runs [\citet{LLC07}].
In practice, if only a single run was made,
we suggest to look at the plot of $\widehat{\bolds{\pi}}$ to choose $k_0$
from where $\widehat{\bolds{\pi}}_k$ has been approximately stable.
Here, we
denote by $\widehat{\bolds{\pi}}_k$
the sampling frequencies of the respective subregions realized by
iteration $k$.
It follows from Theorem \ref{samcconvergence} that $\widehat{\bolds
{\pi}}_k
\rightarrow
\bolds{\pi}$ when the number of iterations, $k$, becomes large.

Trajectory averaging
can directly benefit one's inference in many applications of SAMC. A
typical example is Bayesian
model selection, where the ratio $\omega_i/\omega_j$ just corresponds to
the Bayesian factor of two models if one partitions the sample space
according to the model index and imposes an uniform prior on the model space
as done in \citet{L09}.
Another example is inference for the spatial
models with intractable normalizing constants, for which
\citet{LLC07} has demonstrated how SAMC can be used to
estimate the normalizing constants for these models and how the estimate
can then be used for inference of the model parameters.
An improved estimate of the normalizing constant function would
definitely benefit one's inference for the model.


\section{Trajectory averaging for a stochastic
approximation MLE algorithm}\label{sec4}

Consider the standard missing data problem:
\begin{itemize}
\item$y$ is the observed incomplete data.
\item$f(x,\theta)$ is the complete data likelihood, that is, the
likelihood of
the complete data $(x,y)$ obtained by augmenting the observed data $y$
with the missing data~$x$. The dependence of $f(x,\theta)$ on $y$ is here
implicit.
\item$p(x,\theta)$ is the predictive distribution of the missing data $x$
given the observed data $y$, that is, the predictive likelihood.
\end{itemize}
Our goal is to find the maximum likelihood estimator of $\theta$.
This problem has been considered by a few authors under the framework
of stochastic approximation; see, for example, \citet{Y89}, \citet{GK98}
and \citet{DLM99}.
A basic algorithm proposed by \citet{Y89} for the problem can be
written as
%
%
\begin{equation} \label{MLEeq1}
\theta_{k+1}=\theta_{k}+a_k \partial_{\theta} \log f(X_{k+1},\theta_{k}),
\end{equation}
where the missing data $X_{k+1}$ can be imputed using a MCMC algorithm,
such as
the Metropolis--Hastings algorithm. Under standard regularity
conditions, we have
\[
h(\theta)= E_{\theta}[ \partial_{\theta} \log f(X,\theta) ]
= \partial_{\theta} l(\theta),
\]
where $l(\theta)$ is the log-likelihood function of the
incomplete data.

To show that the trajectory averaging estimator is asymptotically efficient
for a varying truncation version of the algorithm (\ref{MLEeq1}),
we assume (A$_3$), (A$_4$) and some regularity conditions for the
distribution $f(x,\theta)$. The conditions (A$_1$) and (A$_2$) can be
easily verified
with the following settings:

\begin{itemize}
\item The Lyapunov function $v(\theta)$ can be chosen as
$v(\theta)= -l(\theta)+C$,
where $C$ is chosen such that $v(\theta)>0$. Thus,
\[
\langle\nabla v(\theta),h(\theta)\rangle=-\| \partial_{\theta}
l(\theta
) \|^2.
\]
The set of stationary points of (\ref{MLEeq1}), $\{ \theta\dvtx \langle
\nabla v(\theta),h(\theta)
\rangle=0 \}$, coincides with the set of the solutions
$\{\theta\dvtx \partial_{\theta} l(\theta)=0 \}$. Then the condition
(A$_1$) can be verified
by verifying that $l(\theta)$ is continuously differentiable (this is
problem dependent).

\item The matrix $F$ trivially is the Hessian matrix of $l(\theta)$.
Then (A$_2$) can be
verified using the Taylor expansion.
\end{itemize}

In summary, we have the following theorem.
\begin{theorem} \label{MLEefficiency} Assume  conditions \textup{(A$_3$)} and
\textup{(A$_4$)} hold. Then
the estimator $\bar{\theta}_k$ generated by a varying truncation
version of algorithm (\ref{MLEeq1})
is asymptotically efficient.
\end{theorem}

In practice, to ensure the drift condition to be satisfied, we may
follow \citet{AMP05}
to impose some constraints on the tails of the distribution $f(x,\theta
)$ and the proposal
distribution $q(x,y)$. Alternatively, we can follow \citet{LLC07}
to choose a proposal satisfying the local positive condition (\ref
{proposalcons}) and
to restrict the
sample space ${\mathcal X}$ to be compact.
For example, we may set ${\mathcal X}$ to a huge space, say, ${\mathcal X}
=[-10^{100},10^{100}]^{d_x}$.
As a practical matter,
this is equivalent to setting ${\mathcal X}=\mathbb{R}^{d_x}$.

\section{Conclusion} \label{sec5}

In this paper, we have shown that the trajectory averaging estimator is
asymptotically efficient for a general stochastic approximation MCMC
algorithm under mild conditions,
and then applied this result to the stochastic approximation Monte
Carlo algorithm and a stochastic approximation MLE algorithm.

The main difference between this work and the work published in the
literature, for example,
\citet{PJ92} and \citet{C93}, are at the conditions on
the observation noise. In the literature, it is usually assumed
directly that
the observation noise has the decomposition
$\varepsilon_{k}=e_k+\nu_k$, where $\{e_k\}$ forms
a martingale difference sequence and $\nu_k$ is a higher order term of
$o(a_k^{1/2})$.
As shown in Lemma \ref{lem3},
the stochastic approximation MCMC
algorithm does not satisfy this decomposition.

\begin{appendix}
\section{\texorpdfstring{Proofs of Theorems \protect\lowercase{\ref{contheorem}} and
\protect\lowercase{\ref{efftheorem}}}{Proofs of Theorems 2.2 and 2.3}}

Lemma \ref{lem1} is a partial restatement of Proposition
6.1 of \citet{AMP05}.
\begin{lemma} \label{lem1}
Assume  condition \textup{(A$_3$)} holds. Then the following results hold:
\begin{enumerate}[(B$_1$)]
\item[(B$_1$)] For any $\theta\in\Theta$, the Markov kernel
$P_{\theta
}$ has a single
stationary distribution $f_{\theta}$. In addition, $H\dvtx \Theta\times
{\mathcal X}\to\Theta$ is
measurable for all $\theta\in\Theta$, $\int_{{\mathcal X}} \|
H(\theta,x)\|
f_{\theta}(x) \,dx <\infty$.

\item[(B$_2$)] For any $\theta\in\Theta$, the Poisson equation
$u(\theta,x)-P_{\theta}u(\theta,x)
=H(\theta,x)-h(\theta)$ has a solution $u(\theta,x)$, where
$P_{\theta
} u(\theta,x)=
\int_{{\mathcal X}} u(\theta,x') P_{\theta}(x,x') \,d x'$.
There exist a function $V\dvtx {\mathcal X}\to[1,\infty)$ such that $\{x
\in{\mathcal X},
V(x) <\infty\} \ne\varnothing$,
and a constant $\beta\in(0,1]$ such that for any compact subset
$\mathcal{K}
\subset\Theta$, the
following holds:
%
%
\begin{eqnarray} \label{boundeq1}
\mbox{\textup{(i)} } && \sup_{\theta\in\mathcal{K}} \|H(\theta,x)\|_V < \infty,
\nonumber\\
\mbox{\textup{(ii)} }&& \sup_{\theta\in\mathcal{K}}\bigl( \|u(\theta,x)\|_V+\|P_{\theta}
u(\theta
,x)\|_V\bigr) < \infty, \nonumber\\[-8pt]\\[-8pt]
\mbox{\textup{(iii)} } &&\sup_{(\theta,\theta')\in\mathcal{K}\times\mathcal{K}} \|
\theta-\theta
'\|
^{-\beta}
\bigl(\|u(\theta,x)-u(\theta',x)\|_V\nonumber\\
&&\hspace*{78.8pt}\qquad{}+\|P_{\theta} u(\theta,x)-P_{\theta'}
u(\theta',x)\|_V \bigr) <\infty.\nonumber
\end{eqnarray}
\end{enumerate}
\end{lemma}

Lemma \ref{lem1.1} is a restatement of Proposition
5.1 of \citet{AMP05}.
\begin{lemma} \label{lem1.1} Assume  conditions \textup{(A$_1$), (A$_3$)} and
\textup{(A$_4$)} hold.
Let ${\mathcal X}_0 \subset{\mathcal X}$ be such that $\sup_{x \in
{\mathcal X}_0} V(x)$
$<\infty
$ and that
$\mathcal{K}_0 \subset\mathcal{V}_{M_0}$, where $\mathcal{V}_0$ is
defined in \textup{(A$_1$)}.
Then $\sup_{k} E [ V^{\alpha}(X_k)I(k \geq k_{\sigma_s})] < \infty
$, where
$\alpha\geq2$ is defined in condition \textup{(A$_3$)} and $k_{\sigma_s}$ is
defined in
Theorem \ref{lem50}.
\end{lemma}

Lemma \ref{lem2} is a restatement of Corollary 2.1.10 of Duflo (\citeyear{D97}),
pages 46 and~47.
\begin{lemma} \label{lem2} Let $\{S_{ni}, \mathcal{G}_{ni}, 1 \leq i
\leq k_n,
n\geq1 \}$ be a zero-mean,
square-integrable martingale array with differences $\upsilon_{ni}$, where
$\mathcal{G}_{ni}$ denotes the $\sigma$-field. Suppose that the following
assumptions apply:
\begin{longlist}
\item The $\sigma$-fields are nested: $\mathcal{G}_{ni}
\subseteq\mathcal{G}
_{n+1,i}$ for $1 \leq i \leq k_n$,
$n \geq1$.
\item $\sum_{i=1}^{k_n} E(\upsilon_{ni} \upsilon
_{ni}^T|\mathcal{G}
_{n,i-1}) \rightarrow\Lambda$ in probability,
where $\Lambda$ is a positive definite matrix.
\item For any $\varepsilon>0$, $\sum_{i=1}^{k_n} E [ \|\upsilon
_{ni}\|^2
I_{(\|\upsilon_{ni}\| \geq\varepsilon)} | \mathcal{G}_{n,i-1} ]
\rightarrow0$
in probability.
\end{longlist}
Then $S_{n k_n}=\sum_{i=1}^{k_n} \upsilon_{ni} \rightarrow N(0,
\Lambda
)$ in distribution.
\end{lemma}
\begin{df} \label{dfSLLN} For $\varrho\in(0, \infty)$, a sequence
$\{
X_n, n \geq1\}$ of random variables
is said to be residually Ces\`{a}ro $\varrho$-integrable
[$\operatorname{RCI}(\varrho)$, in short] if
\[
\sup_{ n\geq1} \frac{1}{n} \sum_{i=1}^n E |X_i| <
\infty
\]
and
\[
\lim_{n \rightarrow\infty} \frac{1}{n} \sum_{i=1}^n E( |X_i|
-i^{\varrho}) I( |X_i|>i^{\varrho}) =0.
\]
\end{df}

Lemma \ref{lemmaSLLN} is a restatement of Theorem 2.1 of
Chandra and Goswami (\citeyear{CG06}).
\begin{lemma} \label{lemmaSLLN} Let $\{X_n, n \geq1\}$ be a sequence
of nonnegative
random variables satisfying $E(X_i X_j) \leq E(X_i) E(X_j)$ for all $i
\ne j$ and
let $S_n=\sum_{i=1}^n X_i$. If $\{ X_n, n \geq1\}$ is $\operatorname{RCI}(\varrho)$
for some
$\varrho\in(0,1)$, then
\[
\frac{1}{n} [ S_n- E(S_n) ] \rightarrow0\qquad \mbox{in probability}.
\]
\end{lemma}
\begin{lemma} \label{lem3}
Assume  conditions \textup{(A$_1$), (A$_3$)} and \textup{(A$_4$)} hold.
Let ${\mathcal X}_0 \subset{\mathcal X}$ be such that $\sup_{x \in
{\mathcal X}_0} V(x)$
$<\infty
$ and that
$\mathcal{K}_0 \subset\mathcal{V}_{M_0}$, where $\mathcal{V}_0$ is
defined in \textup{(A$_1$)}.
If $k_{\sigma_s}< \infty$, which is defined in Theorem \ref{lem50},
then there exist $\mathbb{R}^d$-valued random processes $\{e_k\}_{k
\geq k_{\sigma_s}}$,
$\{\nu_k\}_{k \geq k_{\sigma_s}}$ and $\{ \varsigma_k \}_{k \geq
k_{\sigma_s}}$ defined
on a
probability space $(\Omega,\mathcal{F},\mathbb{P})$ such that:
\begin{longlist}
\item $\varepsilon_k=e_k+\nu_k+\varsigma_k$ for $k \geq k_{\sigma_s}$.
%
\item $\{e_k\}_{k \geq k_{\sigma_s}}$ is a martingale difference
sequence, and
$\frac{1}{\sqrt{n}} \sum_{k=k_{\sigma_s}}^n e_k \longrightarrow
N(0,Q)$ in
distribution, where
$Q=\lim_{k \rightarrow\infty} E (e_k e_k^T)$.
%
\item $\frac{1}{\sqrt{k}} \sum_{i=k_{\sigma_s}}^k E \| \nu
_i\|
\rightarrow0$, as
$k \rightarrow\infty$.
%
\item $E \|\sum_{i=k_{\sigma_s}}^k a_i \varsigma_i \|
\rightarrow0$, as
$k \rightarrow\infty$.
\end{longlist}
\end{lemma}
\begin{pf}
(i)
Let $\varepsilon_{k_{\sigma_s}}=\nu_{k_{\sigma
_s}}=\varsigma_{k_{\sigma_s}}=0$, and
%
%
\begin{eqnarray} \label{noisedecomeq}
e_{k+1} &=& u(\theta_k,x_{k+1})-P_{\theta_k} u(\theta_k,x_k) ,
\nonumber\\
\nu_{k+1}&=& [ P_{\theta_{k+1}} u(\theta_{k+1},x_{k+1})
-P_{\theta_k} u(\theta_{k},x_{k+1}) ]\nonumber\\
&&{} + \frac{a_{k+2}-a_{k+1}}{a_{k+1}}
P_{\theta_{k+1}} u(\theta_{k+1},x_{k+1}), \\
\tilde{\varsigma}_{k+1}&=&a_{k+1} P_{\theta_k} u(\theta_k,x_k),
\nonumber\\
\varsigma_{k+1}&=&\frac{1}{a_{k+1}} (\tilde{\varsigma}_{k+1}-\tilde
{\varsigma}_{k+2}).\nonumber
\end{eqnarray}
It is easy to verify that (i) holds by noticing the Poisson equation
given in (B$_2$).

{\smallskipamount=0pt
\begin{longlist}[(iii)]
\item[(ii)] By (\ref{noisedecomeq}), we have
\[
E(e_{k+1}|\mathcal{F}_k)=E(u(\theta_k,x_{k+1})|\mathcal
{F}_k)-P_{\theta_k}u(\theta
_k,x_k)=0,
\]
where $\{\mathcal{F}_k\}_{k \geq k_{\sigma_s}}$ is a
family of $\sigma$-algebras satisfying $\sigma\{\theta_{k_{\sigma
_s}}, x_{k_{\sigma
_s}}\}
\subseteq\mathcal{F}_0$ and
$\sigma\{\theta_{k_{\sigma_s}},\theta_{k_{\sigma_s}+1},\ldots
,\theta_k;x_{k_{\sigma
_s}}, x_{k_{\sigma
_s+1}}, \ldots,x_k\} \subseteq
\mathcal{F}_k
\subseteq\mathcal{F}_{k+1}$
for all $k \geq k_{\sigma_s}$. Hence, $\{e_k\}_{k \geq k_{\sigma_s}}$
forms a martingale
difference sequence.

When $k_{\sigma_s}< \infty$,
there exists a compact set $\mathcal{K}$ such that $\theta_k \in
\mathcal{K}$ for all
$k \geq0$.
Following from Lemmas \ref{lem1} and \ref{lem1.1}, $\{e_k \}_{k
\geq k_{\sigma_s}}$ is
$e_k$ is uniformly square integrable with
respect to $k$, and
the martingale $s_n=\sum_{k=1}^n e_k$ is square integrable for all $n$.

By (\ref{noisedecomeq}), we have
%
%
\begin{eqnarray} \label{lem3proofeq1}
E(e_{k+1} e_{k+1}^T|\mathcal{F}_{k})&=&E [ u(\theta_k,x_{k+1})u(\theta
_k,x_{k+1})^T | \mathcal{F}_{k} ]\nonumber \\
&&{}-
P_{\theta_k} u(\theta_k,x_k) P_{\theta_k} u(\theta_k,x_k)^T\\
&\stackrel{\triangle}{=}& l(\theta_k,x_k).\nonumber
\end{eqnarray}

Following from Lemmas \ref{lem1} and \ref{lem1.1}, $\| l(\theta
_k,x_k) \|$ is uniformly integrable
with respect to $k$.
Hence, $\{ l(\theta_k,x_k), k \geq k_{\sigma_s}\}$ is RCI($\varrho$)
for any
$\varrho>0$ (Definition \ref{dfSLLN}).
Since $\{ E(e_{k+1} e_{k+1}^T|\mathcal{F}_{k})-E(e_{k+1} e_{k+1}^T) \}$
forms a martingale difference sequence, the correlation
coefficient
$\mbox{Corr} ( l(\theta_i,x_i)$, $l(\theta_j, x_j) )=0$ for all
$i\ne j$.
By Lem\-ma~\ref{lemmaSLLN}, we have, as $n \rightarrow\infty$,
%
%
\begin{equation} \label{mcmccon}
\frac{1}{n} \sum_{k=k_{\sigma_s}}^n l(\theta_k,x_k) \rightarrow
\frac{1}{n} \sum_{k=k_{\sigma_s}}^n E l(\theta_k,x_k) \qquad\mbox{in
probability}.
\end{equation}

Now we show that $E l(\theta_k,x_k)$ also converges. It follows from
(A$_1$) and (B$_2$) that
$l(\theta,x)$ is continuous in $\theta$. By the
convergence of $\theta_k$, we can conclude that $l(\theta_k,x)$
converges to $l(\theta^*,x)$ for any $x \in{\mathcal X}$.
Following from Lemmas \ref{lem1}, \ref{lem1.1} and
Lebesgue's dominated convergence theorem, $E l(\theta_k,x_k)$
converges to $E l(\theta^*,x)$.
Combining with (\ref{mcmccon}), we obtain
%
%
\begin{equation} \label{mcmccon2}
\frac{1}{n} \sum_{k=k_{\sigma_s}}^n l(\theta_k,x_k) \rightarrow
E l(\theta^*,x)=\lim_{k \rightarrow\infty} E (e_k e_k^T)\qquad \mbox{in
probability}.
\end{equation}
Since $\|e_k\|$ can be uniformly bounded by an integrable function $c
V(x)$, the
Lindeberg condition is satisfied, that is,
\[
\sum_{i=k_{\sigma_s}}^n E \biggl[ \frac{\|e_i\|^2}{n}
I_{({\|e_i\|}/{\sqrt{n}} \geq\varepsilon)} \Big| \mathcal{F}_{i-1} \biggr]
\rightarrow0\qquad \mbox{as $n \rightarrow\infty$}.
\]
Following from Lemma \ref{lem2}, we have $\sum_{i=k_{\sigma_s}}^n
e_i/\sqrt{n}
\rightarrow N(0,Q)$
by identifying $e_i/\sqrt{n}$ to $\upsilon_{ni}$, $n$ to $k_n$, and
$\mathcal{F}_{i}$ to
$\mathcal{G}_{ni}$.

\item[(iii)] By condition (A$_4$), we have
\[
\frac{a_{k+2}-a_{k+1}}{a_{k+1}}= o(a_{k+2}).
\]
By (\ref{noisedecomeq}) and (\ref{boundeq1}), there exists a
constant $c_1$ such that the following inequality holds:
\[
\|\nu_{k+1}\|_V \leq c_1 \|\theta_{k+1}-\theta_{k}\|+ o(a_{k+2})
=c_1 \|a_{k} H(\theta_k, x_{k+1}) \| + o(a_{k+2}),
\]
which implies, by (\ref{boundeq1}),
that there exists a constant $c_2$ such that
%
%
\begin{equation} \label{nunormeq}
\|\nu_{k+1}\|_{V^2} \leq c_2 a_{k}.
\end{equation}
Since $V(x)$ is square integrable, $\nu_{k}$ is uniformly integrable
with respect to $k$ and
there exists a constant $c_3$ such that
\[
\sum_{k=k_{\sigma_s}}^{\infty} \frac{E \|\nu_k \|}{\sqrt{k}} \leq c_3
\sum_{k=k_{\sigma_s}}^{\infty} \frac{a_k}{\sqrt{k}} < \infty,
\]
where the last inequality follows from condition (A$_4$).
Therefore, (iii) holds by Kronecker's lemma.

\item[(iv)]
A straightforward calculation shows that
\[
\sum_{i=k_{\sigma_s}}^k a_i \varsigma_i=-\tilde{\varsigma}_{k+1} =
-a_{k+1} P_{\theta_k} u(\theta_k,x_k).
\]
By Lemmas \ref{lem1} and \ref{lem1.1}, $E \| P_{\theta_k}
u(\theta
_k,x_k) \|$ is uniformly bounded with respect to~$k$.
Therefore, (iv) holds.\qed
\end{longlist}}
\noqed\end{pf}

By Theorem \ref{lem50}, we have
%
%
\begin{equation} \label{itertheta}
\theta_{k+1}-\theta^*=(\theta_{k}-\theta^*)+a_k h(\theta_k)+ a_k
\varepsilon_{k+1}\qquad
\forall k \geq k_{\sigma_s}.
\end{equation}
To facilitate the theoretical analysis for the random process $\{
\theta
_k \}$,
we define a reduced random process:
$\{ \tilde{\theta}_k \}_{k \geq0}$, where
%
%
\begin{equation} \label{redproc}
\tilde{\theta}_k =\cases{
\theta_k+ \tilde{\varsigma}_{k}, &\quad $k > k_{\sigma_s}$, \cr
\theta_k, &\quad $0 \leq k \leq k_{\sigma_s}$,}
\end{equation}
which is equivalent to set $\tilde{\varsigma}_{k}=0$ for all
$k=0,\ldots,k_{\sigma_s}$. For convenience, we also define
%
%
\begin{equation} \label{redprocrev}
\tilde{\varepsilon}_k=e_k+\nu_k, \qquad k > k_{\sigma_s}.
\end{equation}
It is easy to verify that
%
%
\begin{eqnarray} \label{itertheta2}
\tilde{\theta}_{k+1}-\theta^*&=&(I+a_k
F)(\tilde{\theta}_{k}-\theta^*)\nonumber\\[-8pt]\\[-8pt]
&&{}+a_k
\bigl( h(\theta_k)-F(\tilde{\theta}_k-\theta^*) \bigr) +a_k
\tilde{\varepsilon}_{k+1}\qquad
\forall k \geq k_{\sigma_s},\nonumber
\end{eqnarray}
which implies
%
%
\begin{eqnarray} \label{redcedeqiter}
\tilde{\theta}_{k+1}-\theta^*&=&\Phi_{k,k_{\sigma_s}} (\tilde
{\theta}_{k_{\sigma_s}}-\theta^*)+
\sum_{j=k_{\sigma_s}}^k \Phi_{k, j+1} a_j
\tilde{\varepsilon}_{j+1}\nonumber\\[-8pt]\\[-8pt]
&&{}+
\sum_{j=k_{\sigma_s}}^k \Phi_{k,j+1} a_j \bigl( h(\theta_j)-F(\tilde
{\theta}_j-\theta
^*) \bigr)\qquad
\forall k \geq k_{\sigma_s},\nonumber
\end{eqnarray}
where $\Phi_{k,j}=\prod_{i=j}^k (I+a_i F)$ if $k \geq j$, and $\Phi
_{j,j+1}=I$, and
$I$ denotes the identity matrix.

For $\gamma$ specified in (A$_2$) and a deterministic integer
$k_0$, define the stopping time $\mu=\min\{j\dvtx j \geq k_0,
\|\theta_j-\theta^*\| \geq\gamma\}$ if $\|\theta_{k_0} -\theta
^*\| <
\gamma$ and 0
if $\|\theta_{k_0} -\theta^*\| \geq\gamma$. Define
%
%
\begin{equation} \label{setAeq}
A=\{i\dvtx k_{\sigma_s}< k_0 \leq i < \mu\},
\end{equation}
and let $I_A(k)$ denote the indicator
function; $I_A(k)=1$ if $k \in A$ and 0 otherwise.
Therefore, for all $k \geq k_0$,
%
%
\begin{eqnarray} \label{lem6eq1}\quad
&&(\tilde{\theta}_{k+1}-\theta^*) I_A(k+1)
\nonumber\\
&&\qquad=\Phi_{k,k_0}(\tilde{\theta}_{k_0}-\theta^*) I_A(k+1)
+ \Biggl[ \sum_{j=k_0}^k \Phi_{k,j+1} a_j \tilde{\varepsilon}_{j+1} I_A(j) \Biggr] I_A(k+1)
\\
&&\qquad\quad{} + \Biggl[ \sum_{j=k_0}^k \Phi_{k,j+1} a_j \bigl( h(\theta_j)-
F(\tilde{\theta}_j-\theta^*) \bigr) I_A(j) \Biggr] I_A(k+1).\nonumber
\end{eqnarray}
Including the terms $I_A(j)$ in (\ref{lem6eq1}) facilitates our use of
some results published in \citet{C02} in the later proofs, but it does
not change
equality of (\ref{lem6eq1}). Note that if $I_A(k+1)=1$, then
$I_A(j)=1$ for all $j=k_0, \ldots, k$.
%
%
\begin{lemma} \label{lem4}
\textup{(i)} The following estimate takes place:
%
%
\begin{equation} \label{lem4eq0}
\frac{a_j}{a_k} \leq\exp\Biggl( o(1) \sum_{i=j}^k a_i \Biggr)\qquad \forall k \geq j,
\forall j \geq1,
\end{equation}
where $o(1)$ denotes a magnitude that tends to zero as $j \rightarrow
\infty$.
%

{\smallskipamount=0pt
\begin{longlist}[(iii)]
\item[(ii)]
Let $c$ be a positive constant, then there exists another
constant $c_1$ such that
%
%
\begin{equation} \label{lem4eq00}
\sum_{i=1}^k a_i^r \exp\Biggl( -c \sum_{j=i+1}^k a_j \Biggr) \leq c_1\qquad \forall k
\geq1,
\forall r \geq1.
\end{equation}
%
%
\item[(iii)] There exist constants $c_0>0$ and $c>0$ such that
%
%
\begin{equation} \label{lem4eq1}
\|\Phi_{k,j} \| \leq c_0 \exp\Biggl\{-c \sum_{i=j}^k a_i \Biggr\}\qquad \forall k
\geq j,
\forall j \geq0.
\end{equation}

\item[(iv)] Let $G_{k,j}=\sum_{i=j}^k (a_{j-1}-a_i) \Phi_{i-1,j} +
F^{-1} \Phi_{k,j}$. Then
$G_{k,j}$ is uniformly bounded with respect to both $k$ and $j$
for $1 \leq j \leq k$, and
%
%
\begin{equation} \label{lem4eq2}
\frac{1}{k} \sum_{j=1}^k \| G_{k,j} \| \longrightarrow0\qquad \mbox{as $k
\rightarrow\infty$}.
\end{equation}
\end{longlist}}
\end{lemma}
\begin{pf}
Parts (i) and (iv) are a restatement of Lemma 3.4.1 of \citet{C02}.
The proof of part (ii) can be found in the proof of Lemma 3.3.2 of \citet{C02}.
The proof of part (iii) can be found in the proof of Lemma 3.1.1 of
\citet{C02}.
\end{pf}
\begin{lemma} \label{lem6} If  conditions \textup{(A$_1$)--(A$_4$)} hold, then
\[
\frac{1}{a_{k+1}} E \|(\theta_{k+1}-\theta^*) I_A(k+1) \|^2
\]
is uniformly bounded with respect to $k$, where the set $A$ is as
defined in (\ref{setAeq}).
\end{lemma}
%
%
\begin{pf}
By (\ref{redproc}) and (\ref{noisedecomeq}), we have
\begin{eqnarray*}
\frac{1}{a_{k+1}}\|\theta_{k+1}- \theta^*\|^2
&=&\frac{1}{a_{k+1}}\| \tilde{\theta}_{k+1}-\theta^*-\tilde
{\varsigma}_{k+1}\|^2\\
&\leq&
\frac{2}{a_{k+1}}\| \tilde{\theta}_{k+1}-\theta^* \|^2+ 2
a_{k+1} \|
P_{\theta_k} u(\theta_k,x_k) \|^2.
\end{eqnarray*}
Following from (B$_2$) and Lemma \ref{lem1.1},
it is easy to see that $E \|P_{\theta_k} u(\theta_k,x_k) \|^2$ is
uniformly bounded with respect to $k$.
Hence, to prove the lemma,
it suffices to prove that $\frac{1}{a_{k+1}}E \|
(\tilde{\theta}_{k+1}-\theta^*) I_A(k+1) \|^2 $ is uniformly bounded with
respect to $k$.

By (\ref{redproc}), (A$_2$) and (B$_2$), there exist constants
$c_1$ and $c_2$ such that
%
%
\begin{eqnarray} \label{lem6eq2}
&&
\|h(\theta_j)-F(\tilde{\theta}_j-\theta^*)\| I_A(j)\nonumber\\
&&\qquad = \|h(\theta_j)-F(\theta_j-\theta^*)-F \tilde{\varsigma}_j \| I_A(j)
\nonumber\\[-8pt]\\[-8pt]
&&\qquad \leq\|h(\theta_j)-F(\theta_j-\theta^*)\| I_A(j) +
c_2 a_j \|P_{\theta_{j-1}} u(\theta_{j-1}, x_{j-1}) \| \nonumber\\
&&\qquad \leq c_1 \|\theta_j-\theta^*\|^{1+\rho}+
c_2 a_j \|P_{\theta_{j-1}} u(\theta_{j-1}, x_{j-1}) \|.\nonumber
\end{eqnarray}
In addition, we have
%
%
\begin{eqnarray} \label{lem6eq21}
E \|\tilde{\theta}_{k_0}- \theta^*\|^2 I_A(k_0)&=& E \|
\theta_{k_0}-\theta^*+\tilde{\varsigma}_{k_0}\|^2
I_A(k_0)\nonumber\\[-8pt]\\[-8pt]
&\leq& 2 \| \theta_{k_0}-\theta^* \|^2 I_A(k_0) + 2 E \|\tilde
{\varsigma}
_{k_0}\|^2.\nonumber
\end{eqnarray}
It is easy to see from (\ref{boundeq1}) and (\ref{noisedecomeq}) that
$\tilde{\varsigma}_{k_0}$ is square integrable. Hence, following from
(\ref{setAeq}),
there exists a constant $\tilde{\gamma}$ such that
%
%
\begin{equation} \label{eq??1}
E \|\tilde{\theta}_{k_0}- \theta^*\|^2 I_A(k_0) \leq\tilde{\gamma}.
\end{equation}

By (\ref{lem6eq1}), (\ref{lem4eq1}), (\ref{lem6eq2}) and (\ref
{eq??1}), and following
Chen [(\citeyear{C02}), page 141] we have
\begin{eqnarray*}
\hspace*{-5pt}&& \frac{1}{a_{k+1}} E \|(\tilde{\theta}_{k+1}-\theta^*) I_A(k+1)\|^2
\\
\hspace*{-5pt}&&\qquad\leq\frac{5 c_0 \tilde{\gamma}}{a_{k+1}} \exp\Biggl( -2c \sum_{i=k_0}^k
a_i \Biggr) \\
\hspace*{-5pt}&&\qquad\quad{} + \frac{5 c_0^2}{a_{k+1}} \sum_{i=k_0}^k \sum_{j=k_0}^k \Biggl[ \exp
\Biggl(-c\hspace*{-0.5pt}\sum
_{s=j+1}^k a_s
\Biggr) a_j \exp\Biggl(-c\hspace*{-0.5pt}\sum_{s=i+1}^k a_s \Biggr) a_i \|E e_{i+1} e_{j+1}^T\| \Biggr] \\
\hspace*{-5pt}&&\qquad\quad{} + \frac{5 c_0^2}{a_{k+1}} \sum_{i=k_0}^k \sum_{j=k_0}^k \Biggl[ \exp
\Biggl(-c\hspace*{-0.5pt}\sum
_{s=j+1}^k a_s
\Biggr) a_j \exp\Biggl(-c\hspace*{-0.5pt}\sum_{s=i+1}^k a_s \Biggr) a_i E \| \nu_{i+1} \nu_{j+1}^T\|
\Biggr]
\\
\hspace*{-5pt}&&\qquad\quad{} + \frac{5 c_0^2 c_2^2 }{a_{k+1}} \sum_{i=k_0}^k \sum_{j=k_0}^k \Biggl[
\exp
\Biggl(-c\sum_{s=j+1}^k a_s
\Biggr) a_j^2 \exp\Biggl(-c\hspace*{-0.5pt}\sum_{s=i+1}^k a_s \Biggr)\\
\hspace*{-5pt}&&\hspace*{83.3pt}\qquad\quad{}\times a_i^2
E \|P_{\theta_{i-1}} u(\theta_{i-1}, x_{i-1}) (P_{\theta_{j-1}}
u(\theta_{j-1}, x_{j-1}))^T \| \Biggr] \\
\hspace*{-5pt}&&\qquad\quad{} + \frac{5c_0^2 c_1^2}{a_{k+1}}
E \Biggl[ \sum_{j=k_0}^k \exp\Biggl(-c\hspace*{-0.5pt}\sum_{s=j+1}^k a_s \Biggr) a_j
\| \theta_j-\theta^*\|^{1+\rho} I_A(j) \Biggr]^2 \\
\hspace*{-5pt}&&\qquad\stackrel{\triangle}{=} I_1+I_2+I_3+I_4+I_5.
\end{eqnarray*}

By (\ref{lem4eq0}), there exists a constant $c_3$ such that
such that
\[
\|I_1\| \leq\frac{ 5c_0 c_3 \tilde{\gamma}}{a_{k_0}} \exp\Biggl(
o(1)\sum_{i=k_0}^{k+1} a_i \Biggr) \exp\Biggl( -2c \sum_{i=k_0}^k a_i \Biggr),
\]
where $o(1) \rightarrow0$ as $k_0 \rightarrow\infty$. This implies that
$o(1)-2c<0$ if $k_0$ is large enough. Hence, $I_1$ is bounded if $k_0$
is large enough.

By (\ref{lem4eq0}) and (\ref{lem4eq00}), for large enough $k_0$,
there exists a constant $c_4$ such that
%
%
\begin{equation} \label{lem4eq4}
\sum_{j=k_0}^k \frac{a_j^2}{a_{k+1}} \exp\Biggl( -c \sum_{s=j+1}^k a_s \Biggr)
\leq\sum_{j=k_0}^k a_j \exp\Biggl( -\frac{c}{2} \sum_{s=j+1}^k a_s \Biggr)
\leq c_4.
\end{equation}
Since $\{e_i\}$ forms a martingale difference sequence (Lemma \ref{lem3}),
\[
E e_{i} e_j^T =E(E(e_i|\mathcal{F}_{i-1})e_j^T)=0\qquad \forall i>j,
\]
which implies that
\begin{eqnarray*}
I_2&=&\frac{5 c_0^2}{a_{k+1}} \sum_{i=k_0}^k \Biggl[ a_i^2 \exp\Biggl(-2c\sum
_{s=j+1}^k a_s
\Biggr) E\|e_i\|^2 \Biggr] \\
&\leq& 5 c_0^2 \sup_i E\|e_i\|^2
\sum_{i=k_0}^k \Biggl[ a_i^2 \exp\Biggl(-2c\sum_{s=j+1}^k a_s \Biggr) \Biggr].
\end{eqnarray*}
Since $\{\|e_i\|, i \geq1\}$ is uniformly bounded by a function
$cV(x)$ which is square integrable,
$\sup_i E\|e_i\|^2$ is bounded by a constant. Furthermore, by (\ref
{lem4eq00}),
$I_2$ is uniformly
bounded with respect to $k$.

By (\ref{noisedecomeq}), (\ref{boundeq1}) and condition (A$_4$),
there exist a
constant $c_0$ and a constant $\tau\in(0,1)$ such that the following
inequality holds:
%
%
\begin{equation} \label{nunormeq2}\qquad
\|\nu_{k+1}\|_V \leq c_0 \|\theta_{k+1}-\theta_{k}\|+ o(a_{k+2})
\leq c_0 b_k + o(a_{k+2})=O\bigl(a_k^{({1+\tau})/{2}}\bigr).
\end{equation}
This, by (B$_1$) and the Cauchy--Schwarz inequality, further implies
that there exists
a constant $c_0'$ such that
%
%
\begin{equation} \label{nunormeq222}
E \| \nu_{i+1} \nu_{j+1}^T \| \leq c_0' a_i^{ ({1+\tau})/{2}}
a_j^{({1+\tau})/{2}}.
\end{equation}
Therefore, there exists a constant $c_5$ such that
\begin{eqnarray*}
I_3 & = & 5c_0^2 \sum_{i=k_0}^k \sum_{j=k_0}^k \Biggl[ \exp\Biggl(-c \sum
_{s=j+1}^k a_s
\Biggr) \frac{a_j}{\sqrt{a_{k+1}}}\\
&&\hspace*{58pt}{}\times \exp\Biggl(-c \sum_{s=i+1}^k a_s \Biggr)
\frac{a_i}{\sqrt{a_{k+1}}} O\bigl(a_i^{({1+\tau})/{2}}\bigr) O\bigl(a_j^{({1+\tau})/{2}}\bigr) \Biggr] \\
& \leq & 5c_0^2 c_5 \sum_{i=k_0}^k \sum_{j=k_0}^k \Biggl[ \exp\Biggl(-\frac{c}{2}
\sum_{s=j+1}^k a_s
\Biggr) {a_j}^{{1}/{2}}\\
&&\hspace*{67.5pt}{}\times \exp\Biggl(-\frac{c}{2} \sum_{s=i+1}^k a_s \Biggr)
{a_i}^{{1}/{2}} a_i^{({1+\tau})/{2}} a_j^{({1+\tau})/{2}}
\Biggr] \\
& = &5 c_0^2 c_5 \Biggl\{ \sum_{j=k_0}^k \Biggl[
{a_j}^{1+{\tau}/{2}} \exp\Biggl(-\frac{c}{2} \sum_{s=j+1}^k a_s
\Biggr) \Biggr] \Biggr\}^2.
\end{eqnarray*}
By (\ref{lem4eq00}), $I_3$ is uniformly bounded with respect to $k$.

Following from Lemmas \ref{lem1} and \ref{lem1.1},
$ E \| P_{\theta_{i-1}} u(\theta_{i-1}, x_{i-1}) (P_{\theta_{j-1}}
u(\theta_{j-1},\break x_{j-1}))^T \|$
is uniformly bounded with respect to $k$.
Therefore, there exists a constant $c_6$ such that
\begin{eqnarray*}
I_4 & = & 5c_0^2 c_2^2 c_6 \sum_{i=k_0}^k \sum_{j=k_0}^k \Biggl[ \exp\Biggl(-c
\sum
_{s=j+1}^k a_s
\Biggr) \frac{a_j^2}{\sqrt{a_{k+1}}} \exp\Biggl(-c \sum_{s=i+1}^k a_s \Biggr)
\frac{a_i^2}{\sqrt{a_{k+1}}} \Biggr] \\
& \leq & 5 c_0^2 c_2^2 c_6 \Biggl\{ \sum_{j=k_0}^k \Biggl[
{a_j}^{{3}/{2}} \exp\Biggl(-\frac{c}{2} \sum_{s=j+1}^k a_s
\Biggr) \Biggr] \Biggr\}^2.
\end{eqnarray*}
By (\ref{lem4eq00}), $I_4$ is uniformly bounded with respect to $k$.


The proof for the uniform boundedness of $I_5$ can be found in the
proof of Lemma 3.4.3 of
Chen (\citeyear{C02}), pages 143 and 144.
\end{pf}
\begin{lemma} \label{lem7} If  conditions \textup{(A$_1$)--(A$_4$)} hold,
then as $k \rightarrow\infty$,
\[
\frac{1}{\sqrt{k}} \sum_{i=k_{\sigma_s}}^k \| h(\theta_i)
-F(\tilde{\theta}_i-\theta^*) \| \longrightarrow0 \qquad\mbox{in probability}.
\]
\end{lemma}
\begin{pf}
By (\ref{redproc}) and (\ref{noisedecomeq}), there exists a constant
$c$ such that
\begin{eqnarray*}
&&\frac{1}{\sqrt{k}} \sum_{i=k_{\sigma_s}}^k \| h(\theta_i)
-F(\tilde{\theta}_i-\theta^*) \|
\\
&&\qquad\leq \frac{1}{\sqrt{k}} \sum_{i=k_{\sigma_s}}^k \| h(\theta_i)
-F(\theta
_i-\theta^*) \|
+\frac{c}{\sqrt{k}} \sum_{i=k_{\sigma_s}}^k a_i \|P_{\theta_{i-1}}
u(\theta
_{i-1}, x_{i-1}) \|\\
&&\qquad \stackrel{\triangle}{=} I_1+I_2.
\end{eqnarray*}
To prove the lemma, it suffices to prove that $I_1$ and $I_2$ both converge
to zero in probability as $k \rightarrow\infty$.

Following from Lemmas \ref{lem1} and \ref{lem1.1},
$E\|P_{\theta_{k}} u(\theta_{k}, x)\|$ is uniformly bounded for all $k
\geq k_{\sigma_s}$.
This implies, by condition (A$_4$), there exists a constant $c$ such that
\[
\sum_{i=1}^{\infty} \frac{a_i
E \|P_{\theta_{i-1}} u(\theta_{i-1}, x_{i-1})\| }{\sqrt{i}}
< c \sum_{i=1}^{\infty} \frac{a_i}{\sqrt{i}} <\infty.
\]
By Kronecker's lemma, $E(I_2) \rightarrow0$, and thus $I_2 \to0$ in
probability.

The convergence $I_1 \rightarrow0$ can be
established as in Chen [(\citeyear{C02}), Lemma 3.4.4] using the condition (A$_2$)
and Lemma \ref{lem6}.
\end{pf}


\textit{Proof of Theorem \protect\ref{contheorem}.}\quad
By Theorem \ref{lem50}, $\theta_k$ converges to the zero point
$\theta^*$
almost surely and
\[
\theta_{k+1}=\theta_k+a_k H(\theta_k,x_{k+1})\qquad \forall k \geq
k_{\sigma_s}.
\]
Consequently, we have, by (\ref{redproc}),
%
%
\begin{eqnarray} \label{proeq1}
\sqrt{k} (\bar{\theta}_k-\theta^*) 
&=&o(1)+\frac{1}{\sqrt{k}} \sum_{i=k_{\sigma_s}}^k (\theta_i-\theta^*)
\nonumber\\[-8pt]\\[-8pt]
&=&o(1)+\frac{1}{\sqrt{k}} \sum_{i=k_{\sigma_s}}^k (\tilde{\theta
}_i-\theta^*)
-\frac{1}{\sqrt{k}} \sum_{i=k_{\sigma_s}}^k
\tilde{\varsigma}_i,\nonumber
\end{eqnarray}
where $o(1) \rightarrow0$ as $k \rightarrow\infty$.

Condition (A$_4$) implies
$\frac{1}{\sqrt{k}} \sum_{i=k_{\sigma_s}}^k a_i \rightarrow0$ by
Kronecker's lemma.
Following Lemmas \ref{lem1} and \ref{lem1.1}, there exists a
constant $c$ such that
%
%
\begin{equation} \label{proeq2}
\frac{1}{\sqrt{k}} \sum_{i=k_{\sigma_s}}^k E \|\tilde{\varsigma
}_i \|
\leq\frac{c}{\sqrt{k}} \sum_{i=k_{\sigma_s}}^k a_{i+1}
\rightarrow0.
\end{equation}
Therefore, $\frac{1}{\sqrt{k}} \sum_{i=k_{\sigma_s}}^k \tilde
{\varsigma}_i \to0$ in
probability as $k \to\infty$.

By (\ref{redcedeqiter}), (\ref{proeq1}) and (\ref{proeq2}), we have
%
%
\begin{eqnarray} \label{proeq3}\quad
\sqrt{k} (\bar{\theta}_k-\theta^*)& = &o_p(1)+ \frac{1}{\sqrt{k}}
\sum_{i=k_{\sigma_s}}^k \Phi_{i-1,k_{\sigma_s}}
(\tilde{\theta}_{k_{\sigma_s}}-\theta^*) \nonumber\\
&&{} + \frac{1}{\sqrt{k}}
\sum_{i=k_{\sigma_s}}^k
\sum
_{j=k_{\sigma_s}}^{i-1}
\Phi_{i-1,j+1}a_j \tilde{\varepsilon}_{j+1}\nonumber\\[-8pt]\\[-8pt]
&&{} +\frac{1}{\sqrt{k}} \sum_{i=k_{\sigma_s}}^k \sum_{j=k_{\sigma_s}}^{i-1}
\Phi_{i-1,j+1} a_j \bigl(h(\theta_j)-F(\tilde{\theta}_j-\theta^*) \bigr)\nonumber\\
&\stackrel{\triangle}{=}& o_p(1)+ I_1+I_2+I_3,\nonumber
\end{eqnarray}
where $o_p(\cdot)$ means
\[
Y_k=o_p(Z_k)\quad\mbox{if and only if}\quad Y_k/Z_k \to0\qquad\mbox{in probability, as
$k \to\infty$.}
\]

By noticing that $\Phi_{k,j}=\Phi_{k-1,j}+a_k F \Phi_{k-1,j}$, we have
\[
\Phi_{k,j}=I+\sum_{i=j}^k a_i F \Phi_{i-1,j}\quad \mbox{and}\quad
F^{-1} \Phi_{k,j}= F^{-1}+\sum_{i=j}^k a_i \Phi_{i-1,j},
\]
and thus
\[
a_{j-1} \sum_{i=j}^k \Phi_{i-1,j}= \sum_{i=j}^k (a_{j-1}-a_i) \Phi
_{i-1,j} +\sum_{i=j}^k a_i \Phi_{i-1,j}.
\]
By the definition of $G_{k,j}$ given in Lemma \ref{lem4}(iv), we have
%
%
\begin{equation} \label{proeq4}
a_{j-1} \sum_{i=j}^k \Phi_{i-1,j}=-F^{-1}+G_{k,j},
\end{equation}
which implies
\[
I_1=\frac{1}{\sqrt{k} a_{k_{\sigma_s}-1}} (-F^{-1}+G_{k,k_{\sigma_s}})
(\tilde{\theta}_{k_{\sigma_s}}-\theta^*).
\]
By Lemma \ref{lem4}, $G_{k,j}$ is bounded.
Therefore, $I_1 \rightarrow0$ as $k \rightarrow\infty$.
The above arguments also imply that there exists a constant $c_{0}>0$
such that
%
%
\begin{equation} \label{proeq5}
\Biggl\|a_{j} \sum_{i=j+1}^k \Phi_{i-1,j+1} \Biggr\| < c_0\qquad \forall k, \forall j<k.
\end{equation}
By (\ref{proeq5}), we have
\begin{eqnarray*}
\|I_3\| &=& \frac{1}{\sqrt{k}} \Biggl\| \sum_{j=k_{\sigma_s}}^k \sum
_{i=j+1}^k \Phi
_{i-1,j+1} a_j
\bigl( h(\theta_j)-F(\tilde{\theta}_j-\theta^* )\bigr) \Biggr\|
\\
&\leq& \frac{c_0}{\sqrt{k}} \sum_{j=k_{\sigma_s}}^k \| h(\theta
_j)-F(\tilde{\theta}
_j-\theta^*) \|.
\end{eqnarray*}
It then follows from Lemma \ref{lem7} that $I_3$ converges to zero in
probability as
$k \rightarrow\infty$.

Now we consider $I_2$. By (\ref{redprocrev}) and (\ref{proeq4}),
\begin{eqnarray*}
I_2&=&-\frac{F^{-1}}{\sqrt{k}} \sum_{j=k_{\sigma_s}}^k e_{j+1}+\frac
{1}{\sqrt{k}}
\sum_{j=k_{\sigma_s}}^k G_{k,j+1} e_{j+1}\\
&&{} +\frac{1}{\sqrt{k}} \sum_{j=k_{\sigma_s}}^k (-F^{-1}+G_{k,j+1})
\nu_{j+1}\\
&\stackrel{\triangle}{=}& J_1+J_2+J_3.
\end{eqnarray*}
Since $\{e_j\}$ is a martingale difference sequence,
\[
E(e_i^T G_{k,i}^T G_{k,j} e_j)=E[E(e_i|\mathcal{F}_{i-1})^T G_{k,i}^T G_{k,j}
e_j]=0\qquad \forall i>j,
\]
which implies that
\[
E\|J_2\|^2 =\frac{1}{k} \sum_{j=k_{\sigma_s}}^k E ( e_{j+1}^T G_{k,j+1}^T
G_{k,j+1} e_{j+1} )
\leq\frac{1}{k} \sum_{j=k_{\sigma_s}}^k \|G_{k,j+1}\|^2 E \|
e_{j+1}\|^2.
\]
By the uniform boundedness of $\{ E \|e_i\|^2, i \geq k_{\sigma_s}\}$,
(\ref{lem4eq2}) and the uniform boundedness of $G_{k,j}$,
there exists a constant $c_1$ such that
%
%
\begin{equation} \label{proeq6}
E\|J_2\|^2 \leq\frac{c_1}{k} \sum_{j=k_{\sigma_s}}^k \|G_{k,j+1}\|
\rightarrow0\qquad
\mbox{as $k \rightarrow\infty$}.
\end{equation}
Therefore, $J_2 \rightarrow0$ in probability as $k \to\infty$.

Since $G_{k,j}$ is uniformly bounded with respect to both $k$ and $j$,
there exists a constant $c_2$
such that
\[
E \|J_3\| \leq\frac{c_2}{\sqrt{k}} \sum_{j=k_{\sigma_s}}^k E \|\nu
_{j+1}\|.
\]
Following from Lemma \ref{lem3}(iii), $J_3$ converges to zero in
probability as $k \to\infty$.

By Lemma \ref{lem3}, $J_1 \rightarrow N(0,S)$ in distribution.
Combining with
the convergence results of $I_1$, $I_3$, $J_2$ and $J_3$, we conclude
the proof of the theorem.\vspace*{12pt}

\textit{Proof of Theorem \ref{efftheorem}.}\quad
Since the order of $\varsigma_k$ is difficult to treat, we consider the
following
stochastic approximation MCMC algorithm:
%
%
\begin{equation} \label{reveffeq1}
\tilde{\theta}_{k+1}=\tilde{\theta}_{k}+a_k
\bigl( h(\theta_k) +\tilde{\varepsilon}_{k+1} \bigr),
\end{equation}
where $\{\tilde{\theta}_k \}$ and $\{\tilde{\varepsilon}_k \}$ are as
defined in
(\ref{redproc}) and (\ref{redprocrev}), respectively. Following from
Lemma \ref{lem3}(ii),
$\{\tilde{\varepsilon}_k\}$ forms a sequence of asymptotically unbiased
estimator of 0.

Let $\hspace*{1.5pt}\bar{\hspace*{-1.5pt}\tilde{\theta}}_n=\sum_{k=1}^n \tilde{\theta}_k/n$. To
establish that $\hspace*{1.5pt}\bar{\hspace*{-1.5pt}\tilde{\theta}}$
is an asymptotically efficient estimator of $\theta^*$, we will first
show (in step 1)
%
%
\begin{equation} \label{reveffeq4}
\sqrt{n} (\hspace*{1.5pt}\bar{\hspace*{-1.5pt}\tilde{\theta}}-\theta^*) \to N(\mathbf{0},\Gamma),
\end{equation}
where ${\Gamma}=F^{-1}Q (F^{-1})^T$, $F=\partial h(\theta^*)/\partial
\theta$
and $Q=\lim_{k \rightarrow\infty} E (e_k e_k^T)$; and then show (in
step 2)
that the asymptotic covariance matrix of $\sum_{k=1}^n \tilde
{\varepsilon}
_k/\sqrt{n}$
is equal to~$Q$.

 \textit{Step} 1.
By (\ref{redprocrev}), we have
%
%
\begin{equation} \label{Noveq111}
\hspace*{1.5pt}\bar{\hspace*{-1.5pt}\tilde{\theta}}= \bar{\theta} +\frac{1}{n} \sum_{k=1}^n
\tilde
{\varsigma}_k.
\end{equation}
By Lemmas \ref{lem1} and \ref{lem1.1},
$E \|P_{\theta_{k-1}} u(\theta_{k-1},x_{k-1})\| $ is uniformly bounded
for $k \geq k_{\sigma_s}$ and thus
there exists a constant $c$ such that
\[
E \Biggl\| \frac{1}{\sqrt{n}} \sum_{k=k_{\sigma_s}}^n \tilde{\varsigma
}_k \Biggr\|
= E \Biggl\| \frac{1}{ \sqrt{n} } \sum_{k=k_{\sigma_s}}^n a_k P_{\theta_{k-1}}
u(\theta
_{k-1},x_{k-1}) \Biggr\|
\leq\frac{c}{ \sqrt{n}} \sum_{k=k_{\sigma_s}}^n a_k.
\]
By Kronecker's lemma and (A$_4$), we have
$\frac{1}{\sqrt{n}} \sum_{k=k_{\sigma_s}}^n a_k \to0$ in
probability. Hence,
$\frac{1}{\sqrt{n}} \sum_{k=k_{\sigma_s}}^n \tilde{\varsigma}_k$
$=o_p(1)$ and
%
%
\begin{equation} \label{reveffeq2.5}
\frac{1}{n} \sum_{k=k_{\sigma_s}}^n \tilde{\varsigma}_k = o_p(n^{-1/2}).
\end{equation}
That is
%
%
\begin{equation} \label{reveffeq3}
\hspace*{1.5pt}\bar{\hspace*{-1.5pt}\tilde{\theta}}_n
= \bar{\theta}_n +o_p(n^{-1/2}).
\end{equation}
Following from Theorem \ref{contheorem} and Slutsky's theorem, (\ref
{reveffeq4}) holds.

 \textit{Step} 2.
Now we show the asymptotic covariance matrix of $\sum_{k=1}^n
\tilde{\varepsilon}
_k/\sqrt{n}$
is equal to $Q$. Consider
\begin{eqnarray*}
&& E \Biggl( \frac{1}{\sqrt{n}} \sum_{k=1}^n \tilde{\varepsilon}_k \Biggr) \Biggl( \frac
{1}{\sqrt
{n}} \sum_{k=1}^n \tilde{\varepsilon}_k \Biggr)^T -
\frac{1}{n} \Biggl( \sum_{k=1}^n E (\tilde{\varepsilon}_k ) \Biggr) \Biggl( \sum
_{k=1}^n E
(\tilde{\varepsilon}_k ) \Biggr)^T \\
&&\qquad = \frac{1}{n} \sum_{k=1}^n E(\tilde{\varepsilon}_k \tilde{\varepsilon
}_k^T) +
\frac{1}{n} \sum\sum_{i \ne j} E(\tilde{\varepsilon}_i \tilde
{\varepsilon}_j^T)
- \frac{1}{n} \Biggl[ \sum_{k=1}^n E(\tilde{\varepsilon}_k ) \Biggr] \Biggl[ \sum_{k=1}^n
E(\tilde{\varepsilon}_k ) \Biggr]^T \\
&&\qquad = (I_1)+(I_2)+(I_3).
\end{eqnarray*}

By (\ref{redprocrev}), we have
\begin{eqnarray*}
(I_1)&=&
\frac{1}{n} \sum_{k=1}^n E(e_k e_k^T) + \frac{2}{n} \sum_{k=1}^n E(e_k
\nu_k^T) +\frac{1}{n} \sum_{k=1}^n E(\nu_k \nu_k^T)\\
&=&(J_1)+(J_2)+(J_3).
\end{eqnarray*}
By (\ref{nunormeq2}), $ \| \nu_k \nu_k^T \|_{V^2} =O(a_k^{1+\tau})$ for
$k \geq k_{\sigma_s}$, where
$\tau\in(0,1)$ is defined in (A$_4$). Since $V^2(x)$ is square
integrable, there exists a constant $c$ such that
\[
\frac{1}{n} \sum_{k=1}^n E \| \nu_k \nu_k^T \| \leq o(1)+ \frac
{c}{\sqrt
{n}} \frac{1}{\sqrt{n}} \sum_{k=k_{\sigma_s}}^n a_k^{1+\tau},
\]
which, by Kronecker's lemma and (A$_4$), implies $J_3 \to0$ as $n \to
\infty$.

Following from Lemmas \ref{lem1} and \ref{lem1.1}, $\{ \|e_k\| \}
_{k \geq k_{\sigma_s}}$ is uniformly bounded with respect to $k$.
Therefore, there exists a constant $c$ such that
\[
J_2=\frac{2}{n} \sum_{k=1}^n E\|e_k \nu_k^T \| \leq o(1)+ \frac{c}{n}
\sum_{k=k_{\sigma_s}}^n E \| \nu_k\|.
\]
Following from Lemma \ref{lem3}(iii), $J_2 \to0$ as $n \to\infty$.

By (\ref{lem3proofeq1}), $E(e_{k+1} e_{k+1}^T) =E l(\theta_k, x_k)$.
Since $l(\theta,x)$ is continuous in $\theta$, it follows from Theorem
\ref{lem50} that
$l(\theta_k,x)$ converges to $l(\theta^*,x)$ for any $x \in{\mathcal X}$.
Furthermore,
following from Lemma \ref{lem1.1} and Lebesgue's dominated convergence theorem,
we conclude that $E l(\theta_k,x_k)$ converges to $E l(\theta^*,x)$,
and thus
\[
J_1 \to E l(\theta^*,x) =\lim_{k \to\infty} E(e_k e_k^T)=Q.
\]
Summarizing the convergence results of $J_1$, $J_2$ and $J_3$, we
conclude that $(I_1) \to Q$ as $n \to\infty$.

By (\ref{redprocrev}), for $i \ne j$, $i \geq k_{\sigma_s}$ and $j
\geq k_{\sigma_s}$,
we have
%
%
\begin{eqnarray} \label{reveffeq2}
E(\tilde{\varepsilon}_i \tilde{\varepsilon}_j^T) &=&
E \{ ( e_i+ \nu_i )
( e_j+ \nu_j )^T \}=E ( e_i e_j^T+ \nu_i \nu_j^T+e_i \nu_j^T+ \nu_i
e_j^T )\nonumber\\[-8pt]\\[-8pt]
&=& E(\nu_i \nu_j^T),\nonumber
\end{eqnarray}
where the last equality follows from the result that
$\{e_k\}_{k \geq k_{\sigma_s}}$ is a martingale difference sequence
[Lemma \ref{lem3}(ii)].
By (\ref{nunormeq222}), there exists a constant $c$ such that
\[
E \| \nu_i \nu_j^T \| \leq c a_i^{({1+\tau})/{2}} a_j^{({1+\tau})/{2}},
\]
which implies that
%
%
\begin{equation} \label{Noveq1}\qquad \ \
\biggl\| \frac{1}{n} \sum\sum_{i \ne j} E ( \nu_i \nu_j^T) \biggr\| \leq o(1)+c
\Biggl[ \frac{1}{\sqrt{n}} \sum_{i=k_{\sigma_s}}^n a_i^{({1+\tau
})/{2}} \Biggr]
\Biggl[ \frac{1}{\sqrt{n}} \sum_{j=k_{\sigma_s}}^n a_j^{({1+\tau
})/{2}} \Biggr].
\end{equation}
By Kronecker's lemma and (A$_4$), $\sum_{i=k_{\sigma_s}}^n a_i^{(
{1+\tau
})/{2}}/\sqrt{n} \to0$ and thus
%
%
\begin{equation} \label{Noveq2}
\frac{1}{n} \sum\sum_{i \ne j} E ( \nu_i \nu_j^T) \to0\qquad \mbox{as $n
\to\infty$}.
\end{equation}
In summary of (\ref{reveffeq2}) and (\ref{Noveq2}), we have
%
%
\begin{equation} \label{Noveq3}
(I_2)= \frac{1}{n} \sum\sum_{i \ne j} E(\tilde{\varepsilon}_i \tilde
{\varepsilon}_j^T)
\to0\qquad \mbox{as $n \to\infty$}.
\end{equation}

By (\ref{nunormeq2}), there exists a constant $c$ such that
\[
\frac{1}{\sqrt{n}} \Biggl\|\sum_{k=1}^n E \nu_k\Biggr\| \leq o(1)+ \frac
{1}{\sqrt
{n}} \sum_{k=k_{\sigma_s}}^n E \|\nu_k\|
=o(1)+ \frac{c}{\sqrt{n}} \sum_{k=k_{\sigma_s}}^n a_k^{
({1+\tau})/{2}}.
\]
By Kronecker's lemma and (A$_4$), we have
%
%
\begin{equation} \label{Noveq22}
\frac{1}{\sqrt{n}} \Biggl\|\sum_{k=1}^n E \nu_k\Biggr\| \to0\qquad \mbox{as $n \to
\infty$}.
\end{equation}
By Lemma \ref{lem1}(i) and (ii), where it is shown that $\{e_k\}_{k
\geq k_{\sigma_s}}$ is a martingale difference sequence, we have
\begin{eqnarray*}
(I_3)&=& \frac{1}{n} \Biggl[ \sum_{k=1}^n E(e_k +\nu_k) \Biggr] \Biggl[ \sum_{k=1}^n E(e_k
+\nu_k ) \Biggr]^T\\
&=& \Biggl[ \frac{1}{\sqrt{n}} \sum_{k=1}^n E(\nu_k) \Biggr]
\Biggl[ \frac{1}{\sqrt{n}} \sum_{k=1}^n E(\nu_k ) \Biggr]^T.
\end{eqnarray*}
Following from (\ref{Noveq22}), we have $(I_3) \to0$ as $n \to\infty$.

Summarizing the convergence results of $(I_1)$, $(I_2)$ and $(I_3)$,
the asymptotic covariance matrix of $\sum_{k=1}^n \tilde{\varepsilon
}_k/\sqrt{n}$
is equal to $Q$. Combining with (\ref{reveffeq4}), we conclude that
$\hspace*{1.5pt}\bar{\hspace*{-1.5pt}\tilde{\theta}}_k$
is an asymptotically efficient estimator of $\theta^*$.

Since $\hspace*{1.5pt}\bar{\hspace*{-1.5pt}\tilde{\theta}}_k$ and $\bar{\theta}_k$ have the same
asymptotic
distribution $N(\mathbf{0}, \Gamma)$, $\bar{\theta}_k$ is also
asymptotically efficient as an estimator of $\theta^*$.
This concludes the proof of Theorem~\ref{efftheorem}.

\section{\texorpdfstring{Proofs of Theorems \protect\lowercase{\ref{samcconvergence}} and
\protect\lowercase{\ref{samcavetheorem}}}{Proofs of Theorems 3.1 and 3.2}}

The theorems can be proved using Theorems \ref{lem50} and \ref{contheorem}
by showing that SAMC satisfies the conditions
(A$_1$) and (A$_2$), as (A$_3$) is assumed, and (A$_4$) and
and the condition $\sup_{x\in{\mathcal X}_0} V(x)<\infty$ have been verified
in the text.

 \textit{Verification of} (A$_1$).
To simplify notation, in the proof we drop the subscript $k$,
denoting $x_k$ by $x$ and
denote $\theta_k=(\theta_{k}^{(1)},\ldots,\theta_{k}^{(m-1)})$
by $\theta=(\theta^{(1)},\ldots,\theta^{(m-1)})$.
Since the invariant distribution of the MH kernel
is $f_{\theta}(x)$,
we have for any fixed~$\theta$,
%
%
\begin{eqnarray} \label{app22}
E\bigl(I_{\{x\in E_i\}} -\pi_i\bigr)&=&\int_{{\mathcal X}} \bigl( I_{\{x\in E_i\}}-\pi_i\bigr)
f_{\theta}(x) \,d x\nonumber\\
&=&
\frac{ \int_{E_i} \psi(x) \,d x/e^{\theta^{(i)}}}{ \sum_{j=1}^m
[\int_{E_j} \psi(x) \,d x/e^{\theta^{(j)}}]}-\pi_i \\
&=&\frac
{S_i}{S}-\pi_i\nonumber
\end{eqnarray}
for $i=1,\ldots, m-1$,
where $S_i= \int_{E_i} \psi(x) \,d x/e^{\theta^{(i)}}$
and $S=\sum_{i=1}^{m-1} S_i+\int_{E_m} \psi(x) \,dx$.
Therefore,
\[
h(\theta)=\int_{{\mathcal X}} H(\theta,x) f_{\theta}(x) \,dx=
\biggl(\frac{S_1}{S}-\pi_1, \ldots, \frac{S_{m-1}}{S}-\pi_{m-1}\biggr)^T.
\]

It follows from (\ref{app22}) that $h(\theta)$ is a continuous function
of $\theta$.
Let $\Lambda(\theta)=1-\frac{1}{2} \sum_{j=1}^{m-1} (\frac
{S_j}{S}- \pi
_j)^2$, and define
$v(\theta)=-\log(\Lambda(\theta))$ as in (\ref{vfunctioneq}).
As shown below, $v(\theta)$ is continuously differentiable.
Since $ 0 \leq\frac{1}{2} \sum_{j=1}^{m-1}(\frac{S_j}{S}-\pi_j)^2 <
\frac{1}{2} [\sum_{j=1}^{m-1} (\frac{S_j}{S})^2+\pi_j^2)]
\leq1$ for all $\theta\in\Theta$, $v(\theta)$ takes\vspace*{1pt} values in the
interval $[0,\infty)$.

Solving the system of equations formed by (\ref{app22}), we have the
single solution
\[
\theta^{(i)}
=c+\log\biggl(\int_{E_i} \psi(\mathbf{x}) \,d
\mathbf{x}\biggr)-\log(\pi_i),\qquad
i=1, \ldots, m-1,
\]
where $c=-\log(\int_{E_m} \psi(\mathbf{x}) \,d \mathbf{x})+\log(\pi_m)$.
It is obvious that $v(\theta^*)=0$, and $v(\mathcal{L})$ has an empty
interior,
where $\theta^*$ is
specified in Theorem \ref{samcconvergence}. Therefore, (A$_1$)(iv)
is satisfied.

Given the continuity of $v(\theta)$, for any numbers $M_1>M_0>0$,
$\theta^* \in\mbox{int}(\mathcal{V}_{M_0})$, and $\mathcal
{V}_{M_1}$ is a compact set,
where $\mbox{int}(A)$ denotes the interior of the set $A$.
Therefore, (A$_1$)(i) and (A$_1$)(ii) are verified.

To verify the condition (A$_1$)(iii), we have the following calculations:
%
%
\begin{eqnarray} \label{der0}
\frac{\partial S}{\partial\theta^{(i)}}&=& \frac{\partial
S_i}{\partial
\theta^{(i)}} = -S_i, \qquad
\frac{\partial S_i}{\partial\theta^{(j)}} = \frac{\partial
S_j}{\partial\theta^{(i)}}= 0, \nonumber\\[-8pt]\\[-8pt]
\frac{\partial( {S_i}/{S} )}{\partial\theta^{(i)}}&=&-\frac
{S_i}{S}\biggl(1-\frac{S_i}{S}\biggr), \qquad
\frac{\partial( {S_i}/{S} )}{\partial\theta^{(j)}}=
\frac{\partial( {S_j}/{S} )}{\partial\theta^{(j)}}=\frac{S_i
S_j}{S^2}\nonumber
\end{eqnarray}
for $i, j=1, \ldots,m-1$ and $i \ne j $. Let $b=\sum_{j=1}^{m-1}
S_j/S$, then we have
\begin{eqnarray*}
\frac{\partial v(\theta)}{\partial\theta^{(j)}} & = & \frac{1}{2
\Lambda
(\theta)} \sum_{j=1}^{m-1}
\frac{\partial({S_j}/{S}- \pi_j)^2}{\partial\theta^{(j)}} \\
& = &\frac{1}{\Lambda(\theta)} \biggl[
\sum_{j \ne i} \biggl(\frac{S_j}{S}-\pi_j\biggr)
\frac{S_i S_j}{S^2}-\biggl(\frac{S_i}{S}-\pi_i\biggr) \frac{S_i}{S}\biggl(1-\frac
{S_i}{S}\biggr) \biggr]\\
& = & \frac{1}{\Lambda(\theta)} \Biggl[ \sum_{j=1}^{m-1}
\biggl(\frac{S_j}{S}-\pi_j\biggr) \frac{S_i S_j}{S^2}- \biggl(\frac{S_i}{S}-\pi_i\biggr)
\frac
{S_i}{S} \Biggr]
\\
& = &\frac{1}{\Lambda(\theta)} \biggl[
b \mu_{\xi} \frac{S_i}{S}-\biggl(\frac{S_i}{S}-\pi_i\biggr) \frac{S_i}{S} \biggr]
\end{eqnarray*}
for $i=1, \ldots, m-1$, where it is defined $\mu_{\xi}= \sum
_{j=1}^{m-1} (\frac{S_j}{S}-\pi_j)
\frac{S_j}{b S}$. Thus,
%
%
\begin{eqnarray} \label{negder2}
&&\langle\nabla v(\theta), h(\theta) \rangle\nonumber\\
&&\qquad = \frac{1}{\Lambda
(\theta
)} \Biggl[
b^2 \mu_{\xi} \sum_{i=1}^{m-1}
\biggl(\frac{S_i}{S}-\pi_i\biggr) \frac{S_i}{bS} - b
\sum_{i=1}^{m-1} \biggl(\frac{S_i}{S}-\pi_i\biggr)^2 \frac{S_i}{bS} \Biggr] \nonumber
\nonumber\\[-8pt]\\[-8pt]
&&\qquad= -\frac{1}{\Lambda(\theta)} \Biggl[
b \sum_{i=1}^{m-1} \biggl(\frac{S_i}{S}-\pi_i\biggr)^2 \frac{S_i}{bS}- b^2 \mu
_{\xi
}^2 \Biggr] \nonumber\\
&&\qquad= - \frac{1}{\Lambda(\theta)} \bigl(b \sigma_{\xi}^2 +b(1-b) \mu_{\xi
}^2 \bigr)
\leq0,\nonumber
\end{eqnarray}
where $\sigma_{\xi}^2$ denotes the variance of the discrete
distribution defined
in the following table:

\begin{center}
\begin{tabular}{@{}lccc@{}}
\hline
State $(\xi)$ & $\frac{S_1}{S}-\pi_1$ & $\cdots$ & $\frac
{S_{m-1}}{S}-\pi_{m-1}$ \\
[4pt]
Prob. & $\frac{S_1}{bS}$ & $\cdots$ & $\frac{S_{m-1}}{bS}$ \\
\hline
\end{tabular}
\end{center}

If $\theta=\theta^*$, $\langle\nabla v(\theta), h(\theta) \rangle=0$;
otherwise, $\langle\nabla v(\theta), h(\theta) \rangle<0$.
Therefore, (A$_1$)(iii) is satisfied.

 \textit{Verification of} (A$_2$).
To verify this condition, we first show that $h(\theta)$ has bounded
second derivatives.
Continuing the calculation
in (\ref{der0}), we have
\[
\frac{\partial^2 ( {S_i}/{S} )}{\partial(\theta^{(i)})^2 }
=\frac
{S_i}{S} \biggl(1-\frac{S_i}{S}\biggr)
\biggl(1-\frac{2 S_i}{S}\biggr),\qquad
\frac{\partial^2 ( {S_i}/{S} )}{\partial\theta^{(j)}\, \partial
\theta^{(i)}}=
-\frac{S_iS_j}{S^2} \biggl(1-\frac{2S_i}{S}\biggr),
\]
which implies that the second derivative of $h(\theta)$ is uniformly
bounded by noting
the inequality $0< \frac{S_i}{S}<1$.

Let $F=\partial h(\theta)/\partial\theta$. By (\ref{der0}), we have
\[
F=\pmatrix{
-\dfrac{S_1}{S}\biggl(1-\dfrac{S_1}{S}\biggr) & \dfrac{S_1 S_2}{S^2} & \cdots&
\dfrac
{S_1 S_{m-1}}{S^2}
\vspace*{2pt}\cr
\dfrac{S_2S_1}{S^2} & -\dfrac{S_2}{S}\biggl(1-\dfrac{S_2}{S}\biggr) & \cdots&
\dfrac
{S_2S_{m-1}}{S^2}
\vspace*{2pt}\cr
\vdots& \ddots& \vdots& \vdots
\vspace*{2pt}\cr
\dfrac{S_{m-1}S_1}{S^2}& \cdots& \cdots&
-\dfrac{S_{m-1}}{S}\biggl(1-\dfrac{S_{m-1}}{S}\biggr)}.
\]
Thus, for any nonzero vector $\mathbf{z}=(z_1,\ldots,z_{m-1})^T$,
%
%
\begin{eqnarray} \label{SAMCCLTeq2}
\mathbf{z}^T F \mathbf{z}& = & - \Biggl[ \sum_{i=1}^{m-1} z_i^2 \frac
{S_i}{S} - \Biggl(\sum
_{i=1}^{m-1} z_i \frac{S_i}{S}
\Biggr)^2 \Biggr] \nonumber\\
& = & -b \Biggl[ \sum_{i=1}^{m-1} z_i^2 \frac{S_i}{bS} - \Biggl(\sum_{i=1}^{m-1} z_i
\frac{S_i}{bS}
\Biggr)^2 \Biggr]
-b(1-b) \Biggl(\sum_{i=1}^{m-1} z_i \frac{S_i}{bS} \Biggr)^2 \\
& = & -b\operatorname{Var}(Z)-b(1-b) (E(Z) )^2 <0,\nonumber
\end{eqnarray}
where $E(Z)$ and $\operatorname{Var}(Z)$ denote, respectively, the
mean and variance
of the discrete
distribution defined by the following table:

\begin{center}
\begin{tabular}{@{}lccc@{}}
\hline
State $(Z)$ & $z_1$ & $\cdots$ & $z_{m-1}$ \\
[4pt]
Prob. & $\frac{S_1}{bS}$ & $\cdots$ & $\frac{S_{m-1}}{bS}$ \\[1pt]
\hline
\end{tabular}
\end{center}

This implies that the matrix $F$ is negative definite and thus stable.
Applying Taylor's expansion
to $h(\theta)$ at the point $\theta^*$, we have
\[
\| h(\theta)-F (\theta-\theta^*)\| \leq c \|\theta-\theta^* \|
^{1+\rho},
\]
for some constants $\rho\in(0,1]$ and $c>0$, by noting that $h(\theta
^*)=0$ and that
the second derivatives of $h(\theta)$ are uniformly bounded with
respect to $\theta$.
Therefore, (A$_2$) is satisfied.
\end{appendix}

\section*{Acknowledgments}

The author thanks the Editor, Associate Editor and the Referee for
their constructive
comments which have led to significant improvement of this paper.

\printaddresses

\end{document}